\documentclass[11pt]{amsart}
\usepackage{amsmath,amssymb,amsthm,color,tikz}

\newtheorem{theorem}{\sc Theorem}[section]
\newtheorem{lemma}[theorem]{\sc Lemma}
\newtheorem{proposition}[theorem]{\sc Proposition}
\newtheorem{corollary}[theorem]{\sc Corollary}

\theoremstyle{remark}
\newtheorem{remark}[theorem]{Remark}

\newcommand{\be}{\begin{equation}}
\newcommand{\ee}{\end{equation}}
\newcommand{\dlim}{\displaystyle\lim_{n\rightarrow \infty}}

\providecommand{\abs}[1]{\vert#1\vert}
\providecommand{\norm}[1]{\Vert#1\Vert}

\newcommand{\fl}[1]{\lfloor{#1}\rfloor}

\def\cD{\mathcal{D}}

\def\cH{\mathcal{H}}

\def\cW{\mathcal{W}}
\def\cL{\mathcal{L}}

\def\bN{\mathbb{N}}
\def\bP{\mathbb{P}}

\def\bR{\mathbb{R}}
\def\bZ{\mathbb{Z}}

\def\bx{\textbf{x}}
\def\r{\rho}

\def\e{\varepsilon}

 \def\wt{\widetilde}   

\definecolor{darkgreen}{rgb}{0.0,0.5,0.0}
\definecolor{darkblue}{rgb}{0.0,0.0,0.3}
\definecolor{nicosred}{rgb}{0.65,0.1,0.1}
\definecolor{light-gray}{gray}{0.7}
\numberwithin{equation}{section} 

\allowdisplaybreaks[1]

\begin{document}

\title[TASEP with discontinuous jump rates]{Tasep with discontinuous jump rates}
\author[N.~Georgiou]{Nicos Georgiou}
\address{Nicos Georgiou\\ University of Wisconsin--Madison\\ 
Mathematics Department\\ Van Vleck Hall\\ 480 Lincoln Dr.\\  
Madison WI 53706-1388\\ USA.}
\email{georgiou@math.wisc.edu}
\urladdr{http://www.math.wisc.edu/~georgiou}
%\thanks{}
\author[R.~Kumar]{Rohini Kumar}
\address{Rohini Kumar\\ University of California--Santa Barbara\\ 
Statistics Department\\ South Hall 5607A\\   
Santa Barbara CA 93106-3110\\ USA.}
\email{kumar@pstat.ucsb.edu}
\urladdr{http://www.pstat.ucsb.edu/faculty\%20pages/Kumar.htm}

%\urladdr{http://www.math.wisc.edu/~kumar}
\thanks{This project started while Rohini Kumar was a graduate student in 
UW-Madison.}
\author[T.~Sepp\"al\"ainen]{Timo Sepp\"al\"ainen}
\thanks{T.\ Sepp\"al\"ainen was partially supported by 
National Science Foundation grant DMS-0701091 and DMS-1003651, and by the Wisconsin Alumni Research Foundation.} 
\address{Timo Sepp\"al\"ainen\\ University of Wisconsin-Madison\\ 
Mathematics Department\\ Van Vleck Hall\\ 480 Lincoln Dr.\\  
Madison WI 53706-1388\\ USA.}
\email{seppalai@math.wisc.edu}
\urladdr{http://www.math.wisc.edu/~seppalai}
\keywords{Inhomogeneous TASEP, hydrodynamic limit, corner growth model, entropy solution, scalar conservation law, Hamilton-Jacobi equation, discontinuous flux}
\subjclass[2000]{60K35} 
\date{\today}
\begin{abstract}
We prove a hydrodynamic limit for the totally asymmetric simple
exclusion process with spatially inhomogeneous jump rates 
given by a speed function that may admit discontinuities. The limiting
density profiles are described with a variational formula.
This formula enables us to compute explicit density profiles 
even though we have no information about the invariant distributions
of the process.  
In the case of a two-phase flux  for which a suitable 
 p.d.e.\ theory has been developed 
we also observe that the limit profiles are 
entropy solutions of the corresponding scalar conservation law
with a discontinuous speed function. 
\end{abstract}
\maketitle

\section{Introduction}
This paper  studies hydrodynamic limits of totally asymmetric simple
exclusion processes (TASEPs)
 with spatially inhomogeneous jump rates given by functions 
that are allowed to have discontinuities.  We
prove  a general hydrodynamic limit
and compute some explicit solutions, even though information
about invariant distributions is not
available.  The results come through a variational formula that 
 takes 
advantage of the known behavior of the homogeneous TASEP.  This way we
are able to get explicit formulas, even though the usual scenario in
hydrodynamic limits is that explicit equations and solutions require 
explicit computations of expectations under invariant distributions.  
Together with explicit hydrodynamic profiles we can present 
 explicit limit shapes 
for the related last-passage growth models with spatially inhomogeneous 
rates.  

The class of particle  processes we consider are defined by a positive speed function
$c(x)$ defined for $x\in\bR$, lower semicontinuous and assumed to have a
discrete set of discontinuities.  
 Particles reside at sites of $\bZ$, subject 
to the exclusion rule that admits at most one particle at 
each site. 
The dynamical rule is that a particle
jumps from site $i$ to site $i+1$ at rate $c(i/n)$ provided site 
$i+1$ is vacant. Space and time are both scaled by the factor 
$n$ and then we let $n\to\infty$.  We prove  the almost sure vague
convergence of the empirical measure   to a density $\r(x,t)$, assuming that
the initial particle configurations have a well-defined macroscopic density profile
$\r_0$. 

From known behavior of driven conservative particle
systems a natural expectation would be that the 
macroscopic density $\r(x,t)$  of this discontinuous  TASEP ought to
be, in some sense, the unique entropy solution
of an initial value problem  of the type 
\be \r_t + (c(x)f(\r))_x = 0, \quad \r(x,0)= \r_0(x). \label{discpde}\ee
 Our proof of the hydrodynamic limit does
not lead directly to this  scalar conservation law.  
We can make the
connection through some recent PDE theory  in the special case of the two-phase flow
where the speed function is piecewise constant with a single discontinuity.   
In this case   the discontinuous TASEP chooses the unique entropy solution.  
We would naturally expect TASEP to choose the correct entropy solution in general,
but we have not investigated the PDE side of things further   to justify   such a 
  claim.  
 
The remainder of this introduction reviews briefly some relevant literature and then
closes with an overview of the contents of the paper.  The model and the results are
presented in Section \ref{results}. 

\medskip

\textbf{Discontinuous scalar conservation laws.}
The study of  scalar conservation laws  \be \r_t + F(x, \r)_x = 0\label{disc2}\ee
 whose
flux $F$ may admit discontinuities in $x$
has taken off in the last decade.  
As with the multiple weak solutions of even the simplest spatially homogeneous case,
a key issue is the identification of the unique physically relevant solution 
by means of a suitable {\sl entropy condition}.  (See Sect.~3.4 of \cite{evan} for 
textbook theory.)  
Several different entropy 
conditions for the discontinuous case have been proposed, motivated 
by  particular physical problems.  See for example  
\cite{adim-gowd-03, MR2356208, audu-pert-05, chen-even-klin-08, Diehl, Klin-Ris, ostr-02}.  
Adimurthi, Mishra and Gowda 
\cite{MR2356208} discuss   how different theories  
 lead to different choices of relevant solution.  
An interesting phenomenon is that limits of vanishing higher order 
effects can  lead to distinct choices (such as vanishing viscosity vs.\ vanishing
capillarity).   
 
However, the model we study does not offer  more than one choice. 
 In our case the graphs of the different fluxes do not intersect as they are
all  multiples  of  $f(\r)=\r(1-\r)$.  In such  cases   it is expected that all 
the entropy
criteria single out the same solution (Remark 4.4 on p.~811 of
\cite{adim-etal-05}).  By  appeal to the theory developed by 
Adimurthi and Gowda \cite{adim-gowd-03}   we show that the 
discontinuous TASEP chooses entropy solutions of   equation \eqref{discpde}
in the case where $c(x)$ takes two values separated by a 
single discontinuity   

Our approach to the hydrodynamic limit goes via the interface process
whose limit is a Hamilton-Jacobi equation.   
Hamilton-Jacobi equations with discontinuous spatial dependence 
have been studied by 
Ostrov  \cite{ostr-02}  via mollification. 

\medskip

\textbf{Hydrodynamic limits for spatially inhomogeneous, driven  conservative  particle 
systems.}   Hydrodynamic limits for the case where the speed function possesses
some degree of  smoothness 
  were proved over a decade ago by Covert and Rezakhanlou 
 \cite{cove-reza} and   Bahadoran 
  \cite{Bahadoran-98}. For the case where the speed function is continuous, a hydrodynamic limit was proven by Rezakhanlou
in \cite{Reza-2002}  by the method of \cite{sepp99K}.   
  
The most relevant and interesting predecessor to our work is the study of   
  Chen et al.\ \cite{chen-even-klin-08}. They combine an existence proof of
  entropy solutions for \eqref{disc2} under certain technical hypotheses on $F$ 
   with a hydrodynamic limit 
for an attractive zero-range process (ZRP)   with discontinuous 
speed function.   The hydrodynamic limit is proved through a compactness argument for
approximate solutions that utilizes measure-valued solutions.  
The approach follows   \cite{Bahadoran-98, cove-reza} by
 establishing a microscopic entropy inequality which under the limit turns into a
macroscopic entropy inequality.

The scope of \cite{chen-even-klin-08} and  our work are significantly different.  
Our flux $F(x,\r)=c(x)\r(1-\r)$ does not satisfy the hypotheses of \cite{chen-even-klin-08}.  
 Even with spatial inhomogeneities, a ZRP has product-form invariant distributions that can be readily written down
and computed with.    This is a key distinction in comparison with exclusion processes. 
The microscopic entropy inequality in \cite{chen-even-klin-08}
 is derived by a coupling with a stationary process.  

\medskip
  
 Finally, let us emphasize the distinction between the present work and some
 other hydrodynamic limits that feature spatial inhomogeneities.   
Random rates (as for example  in \cite{sepp99K})   
  lead to 
homogenization (averaging) and   the macroscopic   flux does not depend 
on the spatial variable.   Somewhat similar but still fundamentally different is TASEP
with a slow bond.  In this model jumps across bond $(0,1)$ occur at rate $c<1$
while all other jump rates are $1$.  The deep question is whether the slow bond
disturbs the hydrodynamic   profile for all $c<1$.  V.~Beffara, V.~Sidoravicius and M.~E.~Vares 
have announced a resolution of this question in the affirmative.  
% This is known for $c<1/2$ (\cite{jano-lebo-94}, see also Sect.~III.3 in \cite{ligg-99}).
Then the hydrodynamic
limit can be   derived  in the same way as in the main theorem of the present paper.  The solution is not entirely explicit, however: one unknown constant
  remains that quantifies the effect of the slow bond
(see    \cite{sepp01slow}).     \cite{baha-04-aop}   generalizes the hydrodynamic limit of  \cite{sepp01slow}  
  to a broad class of driven particle systems with a microscopic blockage.  
       
\medskip

\textbf{Organization of this paper.}   Section \ref{results} 
contains the main results for the inhomogeneous corner growth model and TASEP.  
 Sections \ref{CGM} and  
\ref{hydrolimit}   prove the   limits.   Section 
\ref{density}  outlines  the explicit computation  of 
  density profiles for the two-phase TASEP.   
Section \ref{pdes}  discusses the  connection with PDE theory.  

\medskip 

\textbf{Notational conventions.}   $\bN = \{1, 2, 3,\dotsc\}$   and   $\bZ_+ = \{0, 1, 2, \dotsc \}$.
The Exp($c$) distribution has 
density $f(x)=ce^{-cx}$ for $0<x<\infty$. 
Two last passage time models appear in our paper:  the   corner growth model 
whose last-passage times are denoted by $G$,  
and the equivalent wedge growth model with last-passage times $T$.  
 $H(x) = \mathbf{1}_{[0,\infty)}(x)$ is the Heavyside function. 
 $C$ is a constant that may change from line to line.

%\medskip

%\textbf{Acknowledgments.} 

\medskip

\section{Results}
\label{results}
The corner growth model connected with TASEP has been a central object of study
in this area since the seminal 1981 paper of Rost \cite{rost}.  So let us begin 
 with an explicit description of the limit shape for a
  two-phase  corner growth model with a boundary along the diagonal.   
Put  independent 
exponential random variables $ \{ Y_v\}_{v\in\bN^2}$ on the points of the lattice 
with distributions 
\be
Y_{(i,j)}\sim\left\{
\begin{array}{lll}

\vspace{0.1 in}
\textrm{Exp}(c_1), & \textrm {if } i < j \\

\vspace{0.1 in}
\textrm{Exp}(c_1\wedge c_2), & \textrm {if } i=j \\

\textrm{Exp}(c_2), &\textrm {if } i > j.
\end {array}
\right.
\ee 
We  assume that the rates satisfy $c_1 \geq c_2.$  

Define the last passage time 
\be G(m,n) = \max_{\pi \in \Pi(m,n)}\sum_{v\in \pi}Y_v, \quad (m,n)\in \bN^2, \label{basicref}\ee 
where $\Pi(m,n)$ is the collection of  weakly increasing nearest-neighbor
 paths in the rectangle $[m]\times[n]$ that start from $(1,1)$ and go up to $(m,n).$
 That is, elements of $\Pi(m,n)$ are sequences 
$\{(1,1)=v_1, v_2, \dotsc, v_{m+n-1}=(m,n)\}$ such that 
$v_{i+1}-v_i=(1,0)$ or $(0,1)$.   
 
\begin{theorem}
\label{two-phaselastpassagetime}
Let the rates $c_1 \geq c_2>0$.  Define  $c =  {c_1}/{c_2}\geq 1$
 and $b  = 2c-1-2\sqrt{c(c-1)}.$ Then the a.s.\ limit 
\[ \Phi (x,y) =\dlim n^{-1}G (\fl{nx},\fl{ny})\] exists
for  all $(x,y)\in(0,\infty)^2$  and is given by
\[
\Phi (x,y) = 
 \begin{cases} c_1^{-1}\left( \sqrt{x}+\sqrt{y}\right)^2 , &\text{if } 0 < x \leq b^2y \\[11pt]
x\dfrac{4c-(1+b)^2}{c_1(1-b^2)}   +y\dfrac{(1+b)^2-4cb^2}{c_1(1-b^2)} , 
&\text{if } b^2y< x < y \\[11pt]
c_2^{-1}\left( \sqrt{x}+\sqrt{y}\right)^2 , & \text{if } y \leq x < +\infty. 
\end{cases}
 \] 
\end{theorem} 

This theorem will be obtained as a side result of the development in
Section \ref{CGM}. 

%\begin{remark}
%Macroscopically, the maximal path can be one of two kinds; a straight line in the homogenous $c_1$ or rate $c_2$ region or a piecewise linear path starting from $(0,0)$ up to $(u,u)$, $u = (x-b^2y)(1-b^2)^{-1}$ and then a staight line from $(u,u)$ to $(x,y).$ The result of this optimization is the second branch of the limiting last passage time constant in Theorem \ref{two-phaselastpassagetime}
%\end{remark}

\begin{center}
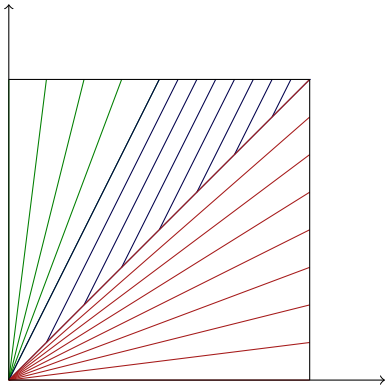
\begin{figure}
\begin{tikzpicture}
\foreach \x in {0,...,4}
\draw[darkgreen] (0,0)--(0.5*\x,4); 

\foreach \x in {0,...,8}
\draw[darkblue] (0.5*\x,0.5*\x)--(2+0.25*\x,4); 

\draw[darkblue] (0,0)--(4,4);

\foreach \x in {0,...,8}
\draw[nicosred] (0,0)--(4,4 - 0.5*\x); 

\draw[->](0,0)--(5,0);
\draw[->](0,0)--(0,5);
\draw(0,4)--(4,4);
\draw(4,0)--(4,4);
\end{tikzpicture}
\caption{Optimal macroscopic paths that give the last 
passage time constant described in Theorem 
\ref{two-phaselastpassagetime}.} 
\end{figure}
\end{center}

We turn to the general hydrodynamic limit.  The variational description needs the
following ingredients.  Define the wedge 
 $$\cW = \{ (x,y)\in \bR^2: y\geq0, \, x\geq -y \}$$ 
 and on $\cW$  the last-passage  function of homogeneous TASEP by 
  \be\gamma(x,y) =(\sqrt{x+y} +\sqrt{y})^2  .\label{gammanot}\ee 
Let $\textbf{x}(s) = (x_1(s), x_2(s))$ denote a path in $\bR^2$ and set 
\begin{align}
\cH(x,y) =\{ \textbf{x}\in C( [0,1], \cW) &: \textbf{x} \text{ is   piecewise }C^1, 
\textbf{x}(0)=(0,0),\notag \\ 
             &\textbf{x}(1)=(x,y), \, \textbf{x}'(s)\in \cW \textrm{ wherever the derivative is defined} \}. \notag
\end{align}  

The  {\sl speed function} $c$ of our system is by assumption 
  a positive lower semicontinuous function on $\bR$. We assume that at each $x\in \bR$ 
\be 
c(x) = \min\Big\{ \lim_{y \nearrow x}c(y),  \lim_{y \searrow x}c(y)  \Big\}.
\label{llrl}
\ee  In particular we assume that the limits in  \eqref{llrl}  exist. We also assume that $c(x)$  has only 
   finitely many discontinuities in any compact set, hence  it
  is bounded away 
  from $0$ in any compact set.   
  
For the hydrodynamic limit   consider
 a sequence of exclusion processes $\eta^n = (\eta^n_i(t): i \in \bZ, \, t\in\bR_+)$ 
  indexed by $n \in \bN$.
These processes
are constructed on a common probability 
space that supports the initial configurations $\{\eta^n(0)\}$ 
and the Poisson clocks of each process.  As always, the clocks of
process $\eta^n$ are independent of its initial state $\eta^n(0)$.
The joint distributions across the index $n$ are immaterial,
except for the assumed initial law of large numbers \eqref{weaklaw}
below.       In the process $\eta^n$ 
a particle at site $i$ attempts  
 a  jump  to $i+1$ with rate $c(i/n)$. Thus the generator of   $\eta^n$ is 
\be
L_nf(\eta) = \sum_{x\in \bZ} c(xn^{-1})\eta(x)(1-\eta(x+1))(f(\eta^{x,x+1})-f(\eta))  
\ee 
for  cylinder functions $f$ on the state space $\{0,1\}^\bZ$.  The usual notation is 
that particle configurations are denoted by $\eta=(\eta(i):i\in\bZ)\in\{0,1\}^\bZ$ and 
\be \eta^{x,x+1}(i) =
\begin{cases}
0 & \text{when } i=x\\
1 & \text{when } i=x+1 \\
\eta(i) &\text{when } i\neq x,x+1 
\end{cases}
\notag
\ee
is the configuration that results from moving a particle from $x$ to $x+1$. 
Let $J_i^n(t)$ denote the number of particles that have made the jump from site $i$ to site $i+1$ in  time interval $[0,t]$ in the process
 $\eta^n.$

An initial macroscopic profile $\rho_0$ is a   measurable function on $\bR$ such that $0\leq \rho_0(x)\leq 1$ for all real $x$, with antiderivative $v_0$ satisfying
\be
v_0(0)=0, \quad v_0(b) - v_0(a) = \int_{a}^{b}\rho_0(x)\,dx.
\label{vdef}
\ee

The macroscopic flux function of the constant rate 1 TASEP  is   \be
f(\r) =\left\{
\begin{array}{ll}

\vspace{0.1 in}
\r(1-\r), & \textrm {if } 0 \leq \r \leq 1 \\ 

-\infty, & \textrm{otherwise. }
\end {array}
\right.
\label{fnot}
\ee 
Its Legendre conjugate   $$f^*(y) = \inf_{0\leq\rho\leq1}\{y\r - f(\r)\}$$ 
represents the limit shape in the wedge. 
We orient our model so that
 growth in the wedge proceeds upward, and so 
we use  $g(y) = -f^*(y)$. It is explicitly given by
\be
g(y) = \sup_{0\leq\rho\leq1}\{f(\rho) - y\rho\} = \left\{
\begin{array}{lll}

\vspace{0.1 in}
-y, & \textrm {if } y \leq -1  \\

\vspace{0.1 in}
\frac1{4}{(1-y)^2}, & \textrm {if } -1 \leq y \leq 1 \\ 

0, & \textrm{if } y \geq 1.
\end {array}
\right.
\ee

For $x \in \bR$  define $v(x,0) = v_0(x)$, and for $t>0,$
\be
v(x,t)=\sup_{w(\cdot)}\left\{  v_0(w(0)) - \int_{0}^{t}c(w(s))\,g\left( \dfrac{w'(s)}{c(w(s))} \right)\,ds \right\}
\label{velocityversion1}
\ee 
where the supremum is taken over continuous piecewise $C^1$ paths $w:[0,t]\longrightarrow \bR$ that satisfy $w(t)=x.$ 
The function $v(x,t)$ is Lipschitz continuous jointly in $(x,t).$
 (see Section \ref{hydrolimit}) and it has a derivative almost everywhere. The macroscopic density is defined by  $\r(x,t) = v_x(x,t).$

The initial distributions of the processes
$\eta^n$  are arbitrary subject to the condition that the following strong law of large numbers holds
at time $t=0$: 
for all real $a<b$
\be
\dlim \frac{1}{n}\sum_{i=\fl{na}+1}^{\fl{nb}}\eta^n_i(0)=\int_a^b\rho_0(x)\,dx \quad \text{a.s.}  
\label{weaklaw}
\ee  

The second theorem gives the hydrodynamic limit of current and particle density for TASEP with discontinuous jump rates.

\begin{theorem}
Let $c(x)$ be a lower semicontinuous positive function satisfying \eqref{llrl}, with finitely many discontinuities in any compact set. Under assumption \eqref{weaklaw}, these strong laws of large numbers hold at each   $t>0$: for all real numbers $a<b$
\be
\dlim n^{-1}J^n_{\fl{na}}(nt) = v_0(a)-v(a,t)\quad a.s.
\label{macrocurrent}
\ee 
and 
\be
\dlim \frac{1}{n}\sum_{i=\fl{na}+1}^{\fl{nb}}\eta^n_i(nt)=\int_a^b\rho(x,t)\,dx \quad a.s.  
\label{macrospeed}
\ee 
where $v(x,t)$ is defined by \eqref{velocityversion1} and $\rho(x,t)=v_x(x,t).$
\label{hydro}
\end{theorem}

\begin{remark} In a totally asymmetric
$K$-exclusion with speed function $c$ the state space would be $\{0,1,\dotsc,K\}^\bZ$
with $K$ particles   allowed at each site, and  one particle   moved from site $x$
to $x+1$ at rate $c(x/n)$ whenever such a move can be legitimately completed.  
 Theorem \ref{hydro} can   be proved for this process with the same method of
 proof. The definition of the limit  \eqref{velocityversion1} would be the same, 
 except that   the explicit  flux $f$ and wedge shape $g$ would be replaced by the
 unknown functions $f$ and $g$ whose existence was proved in \cite{sepp99K}. 
\end{remark}

 To illustrate Theorem \ref{hydro} we compute the macroscopic 
density profiles $\r(x,t)$  
  from constant   initial conditions
in the  two-phase  model with
speed function 
\be     c(x) = c_1(1-H(x))+
c_2H(x) \label{defca}\ee 
where $H(x) = \mathbf{1}_{[0,\infty)}(x)$ is the Heavyside function
and  $c_1 \geq c_2$.  (The case $c_1<c_2$ can then be deduced from 
 particle-hole duality.)  The particles hit 
the region of lower speed as they pass the origin from left to right.  Depending on the
initial density $\rho$, we see the system adjust to this 
discontinuity  in  different ways to match the actual throughput
of particles on either side of the origin.   The maximal
flux on the right is $c_2/4$ which is realized on the left at 
densities $\r^*$ and $1-\r^*$ with 
\[\r^* = \tfrac 12 - \tfrac 12 \sqrt{1-c_2/c_1}.\]
 
\begin{corollary} Let $c_1 \geq c_2$ and
the speed function as in \eqref{defca}.   Then the macroscopic density profiles with initial conditions $\r_0(x,0)=\r$ are given as follows. 

{\rm (i)} Suppose $0< \r < \r^*$. Define $r^*=r^*(\r)= \tfrac 12 - \tfrac 12 \sqrt{1-4\r(1-\r)c_1/c_2} $.  Then 
\be
\r(x,t) = 
\left\{
\begin{array}{llll}

\vspace{0.1 in}
\r & \textrm {if } -\infty \leq x \leq 0\\

\vspace{0.1 in}
r^*  &\textrm{if } 0\leq x \leq c_2(1-2r^*)t \\

\vspace{0.1 in}
\dfrac{1}{2}\left(1- \dfrac{x}{tc_2}\right) &\textrm{if } c_2(1-2r^*)t\leq x \leq c_2(1-2\r)t\\

\r &\textrm{if } (1-2\r)tc_2\leq x < +\infty
\end{array}
\right.
\ee

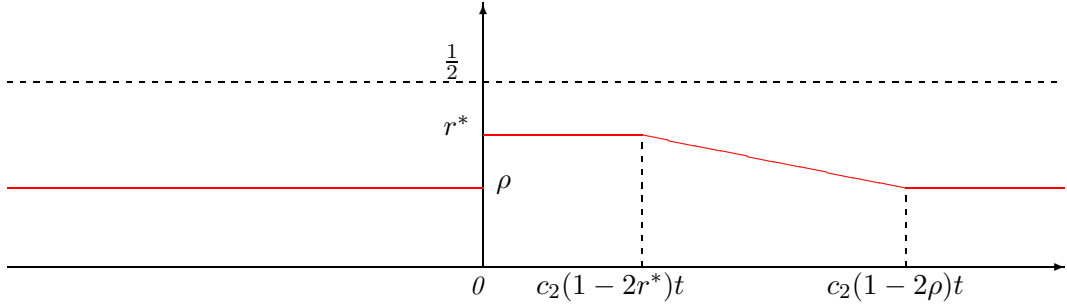
\begin{figure}[ht]
\begin{center}
\begin{picture}(400,130)(20,0)
\put(20,20){\vector(1,0){400}} %axis in the direction of (1,0), length 160 pixels 
\put(200,20){\vector(0,1){100}} %axis in the direction of (0,1) 

{\color{red}
\put(20,50){\line(1,0){180}}
\put(200,70){\line(1,0){60}}
\put(260,70){\line(5,-1){100}}
\put(360,50){\line(1,0){60}}
}
\put(185,70){$r^*$}
\put(205,50){$\r$}
\put(185,95){$\tfrac 12$}
%Dashed lines
\multiput(260,20)(0,5){10}{\line(0,1){2}}
\multiput(360,20)(0,5){6}{\line(0,1){2}}
\multiput(20 ,90)(5,0){80}{\line(1,0){2}}
\put(220,10){$c_2(1-2r^*)t$}
\put(330,10){$c_2(1-2\r)t$}
%Coordinate system x-axis
\put(195,10){\small 0} 
%Coordinate system y-axis

%\put(165,110){$\r(x,t)$}

\end{picture}
\end{center}
\caption{Density profile $\r(x,t)$ in the two-phase ($c_1 > c_2$) TASEP when we start from constant initial configurations $\r_0(x)\equiv \r \in (0,\r^*).$}
\label{fig1}\end{figure}

{\rm (ii)} Suppose $\r^* \leq \r\leq \frac{1}{2}. $ Then 
\be
\r(x,t) = 
\left\{
\begin{array}{llll}

\vspace{0.1 in}
\r & \textrm {if } -\infty \leq x \leq -tc_1(\r-\r^*) \\

\vspace{0.1 in}
1-\r^*  &\textrm{if } -tc_1(\r-\r^*)\leq x \leq 0 \\

\vspace{0.1 in}
\dfrac{1}{2}\left(1-\dfrac{x}{tc_2}\right) &\textrm{if } 0 \leq x \leq (1-2\r)tc_2\\

\r &\textrm{if } (1-2\r)tc_2\leq x < +\infty
\end{array}
\right.
\ee

\begin{figure}[ht]
\begin{center}
\begin{picture}(400,130)(20,0)
\put(20,20){\vector(1,0){400}} %axis in the direction of (1,0), length 160 pixels 
\put(200,20){\vector(0,1){100}} %axis in the direction of (0,1) 

{\color{red}
\put(20,50){\line(1,0){120}}
\put(140,100){\line(1,0){60}}
\put(200,70){\line(5,-1){100}}
\put(300,50){\line(1,0){120}}
}
\put(205,100){$1-\r^*$}
\put(205,50){$\r$}
\put(185,75){$\tfrac 12$}
%Dashed lines
\multiput(140,20)(0,5){16}{\line(0,1){2}}
\multiput(300,20)(0,5){6}{\line(0,1){2}}
\multiput(20 ,70)(5,0){80}{\line(1,0){2}}
\put(110,10){$c_1(\r^*-\r)t$}
\put(260,10){$c_2(1-2\r)t$}
%Coordinate system x-axis
\put(195,10){\small 0} 
%Coordinate system y-axis
%\put(165,110){$\r(x,t)$}

\end{picture}
\end{center}
\caption{Density profile $\r(x,t)$ in the two-phase ($c_1 > c_2$) TASEP when we start from constant initial configurations $\r_0(x)\equiv \r \in [\r^*, \tfrac 12].$}
\label{fig2}\end{figure}

{\rm (iii)} Suppose $ \r\ge \frac{1}{2}.$  Define $r^*=r^*(\r)= \tfrac 12 - \tfrac 12 \sqrt{1-4\r(1-\r)c_2/c_1} $.
Then
\be
\r(x,t) = 
\left\{
\begin{array}{lll}

\vspace{0.1 in}
\r & \textrm {if } -\infty \leq x \leq -tc_1(\r-r^*) \\

\vspace{0.1 in}
1-r^*  &\textrm{if } -tc_1(\r-r^*)\leq x \leq 0 \\

\r &\textrm{if } 0<x <+\infty
\end{array}
\right.
\ee
\label{densityprofiles}

\begin{figure}[ht]
\begin{center}
\begin{picture}(400,130)(20,0)
\put(20,20){\vector(1,0){400}} %axis in the direction of (1,0), length 160 pixels 
\put(200,20){\vector(0,1){100}} %axis in the direction of (0,1) 

{\color{red}
\put(20,70){\line(1,0){120}}
\put(140,100){\line(1,0){60}}
\put(200,70){\line(1,0){220}}
}
\put(205,100){$1-r^*$}
\put(190,65){$\r$}
\put(190,40){$\tfrac 12$}
%Dashed lines
\multiput(140,20)(0,5){16}{\line(0,1){2}}
\multiput(20 ,50)(5,0){80}{\line(1,0){2}}
\put(110,10){$c_1(r^*-\r)t$}

%Coordinate system x-axis
\put(195,10){\small 0} 
%Coordinate system y-axis
%\put(165,110){$\r(x,t)$}

\end{picture}
\end{center}
\caption{Density profile $\r(x,t)$ in the two-phase ($c_1 > c_2$) TASEP when we start from constant initial configurations $\r_0(x)\equiv \r \in (\tfrac 12,1).$}
\label{fig3}\end{figure}
\end{corollary} 

\begin{remark}
Taking $t\to\infty$ in the three cases of Corollary 
\ref{densityprofiles} gives a family of macroscopic profiles that are
fixed by the time evolution.  A natural question to investigate would be the
existence and uniqueness of invariant distributions that correspond
to these macroscopic profiles.  
\end{remark}

Next we relate the density profiles picked by the discontinuous TASEP
to entropy conditions for  scalar conservation laws with discontinuous 
fluxes.
The entropy conditions 
defined by Adimurthi and Gowda \cite{adim-gowd-03} are particularly suited to
our needs. Their results give   
uniqueness of the solution for the scalar conservation law 
\be
\left\{
\begin{array}{ll}

\vspace{0.1 in}
\r_t+\left( F(x,\r) \right)_x = 0, &  x\in \bR, t > 0  \\

\r(x,0) = \r_0(x), &  x \in  \bR
\end {array}
\right.
\label{scl}
\ee with distinct fluxes on the half-lines:
\be
F(x,\r) = H(x)f_r(\r)+(1-H(x))f_{\ell}(\r)  \label{Disflux}
\ee 
where  $f_r,f_{\ell} \in C^{1}(\bR)$ are strictly concave with superlinear decay to $-\infty$
as $\abs{x}\to\infty$.
A solution of \eqref{scl}   means a weak solution, that is,  
$\r \in L^{\infty}_{\text{loc}}(\bR\times\bR_+)$ such that for all continuously differentiable, compactly supported test functions $\phi \in C_c^1\left(\bR\times\bR_+\right)$, 
\be
\int_{-\infty}^{+\infty}\int_{0}^{+\infty}\left( \r\frac{\partial\phi}{\partial t}+F(x,\r)\frac{\partial\phi}{\partial x} \right)\,dt\,dx + \int_{-\infty}^{+\infty}\r(x,0)\phi(x,0)\,dx=0.
\label{weakformulation}
\ee  
\eqref{weakformulation} 
is the weak formulation of the problem
\be
 \begin{cases} 
 \r_t+f_r(\r)_x =0, & \textrm {for } x > 0, t>0  \\[5pt]
 \r_t+f_{\ell}(\r)_x =0, & \textrm{for } x <0, t>0 \\[5pt]
f_r(\r(0+,t)) = f_{\ell}(\r(0-,t))   &\text{for a.e.\ $t>0$} \\[5pt] 
\r(x,0)=\r_0(x). 
\end{cases}
 \label{Breakscl}
\ee 
%  Rankine-Hugoniot (RH) condition   at $x=0$ it satisfies the 
 
The   entropy conditions used in   \cite{adim-gowd-03} come in two sets and assume
the existence of certain one-sided limits: 

\medskip

\textsl{$(E_i)$ Interior entropy condition, or Lax-Oleinik entropy condition:}  
\be
  \r(x+,t) \geq  \r(x-,t) \quad \text{for $x\ne0$ and for all $t>0$.}  
\label{Ei}
\ee  

\textsl{$(E_b)$  Boundary entropy condition at $x=0$:}
for almost every $t$,    the limits $\r(0\pm ,t)$ exist  and  one of the following   holds: 
\be
f_r'(\r(0+,t))\geq 0 \quad\textrm{and}\quad  f_{\ell}'(\r(0-,t))\geq 0,
\label{Eb1}
\ee 
\be
f_r'(\r(0+,t))\leq 0 \quad\textrm{and}\quad  f_{\ell}'(\r(0-,t))\leq 0,
\label{Eb2}
\ee 
\be
f_r'(\r(0+,t))\leq 0 \quad\textrm{and}\quad  f_{\ell}'(\r(0-,t))\geq 0.
\label{Eb3}
\ee 

\smallskip

Define 
\[G_x(p)= \mathbf{1}\{x> 0\}f_r^*(p)  +\mathbf{1}\{x< 0\} f_{\ell}^*(p) + 
\mathbf 1\{x=0 \} \min\bigl( f^*_r(0), f^*_{\ell}(0)\bigr),\]
 where $f_r^*$ and  $f_{\ell}^*$ are the convex duals of $f_r$ and $f_{\ell}$.
Set $V_0(x) = \int_0^x \r_0(\theta)\,d\theta$ and define
\be
V(x,t) = \sup_{w(\cdot)}\left\{ V_0(w(0)) + \int_0^t G_{w(s)}\bigl( w'(s) \bigr)\,ds \right\}
\label{defV} \ee
where the supremum is   over continuous, piecewise linear paths $w: [0,t] \longrightarrow \bR$ with $w(t) = x.$  
 
 \begin{theorem} \cite{adim-gowd-03}
Let $\r_0\in L^{\infty}(\bR)$ and define $V$ by \eqref{defV}.
 Then   $V$ is a uniformly Lipschitz continuous function and 
 $\r(x,t)=V_x(x,t)$ is the unique weak solution of  \eqref{Breakscl} 
 that satisfies the entropy assumptions $(E_i)$ and  
 $(E_b)$ in the class $L^{\infty}\cap BV_{\text{\rm loc}}$ and with discontinuities 
 given by a discrete set of Lipschitz curves. 
 \label{Uniq}
 \end{theorem}

It is easy to check that the two-phase density profile $\r(x,t)$ in Corollary \ref{densityprofiles} is a weak solution (in the sense of \eqref{weakformulation}) to the scalar conservation law \eqref{scl} with flux function $F(x,\r) = c(x)\r(1-\r)$. However we cannot immediately apply this theorem in our case since the two-phase flux function $\widetilde{F}(x,\r) =(1- H(x))c_1f(\r) +  H(x)c_2f(\r)$ is finite only for $\r \in [0,1]$ and in particular is not $C^1.$ We show how we can replace $F(x,\r)$ with $\widetilde{F}(x,\r)$
in the above theorems in Section \ref{pdes}. In particular, we prove the following.

\begin{theorem}
For $\r \in \bR$ define $f_r(\r) = c_2(1-\r)\r$ and $f_{\ell}(\r) = c_1(1-\r)\r$ to be the flux functions for the scalar conservation law \eqref{Breakscl}. 
Let  the initial macroscopic profile 
for the hydrodynamic limit be a measurable function
 $0\le\r_0(x)\le 1$. 
  Then the  macroscopic density profile $\r(x,t)$ from the 
hydrodynamic limit in Theorem \ref{hydro} is the unique
solution described in 
 Theorem \ref{Uniq}.  
\label{pdecor}
\end{theorem} 

\medskip

\section{Wedge last passage time}
\label{CGM}
The strategy of the proof of the hydrodynamic limit is the one from \cite{sepp99K}
and \cite{sepp01slow}.
Instead of the particle process we work with the height process.  
The limit is first proved for the jam initial condition of TASEP (also called step initial
condition) which for the height process is an initial wedge shape.  
This process can be equivalently represented by the wedge  last-passage model.
 Subadditivity 
gives the limit.   The general case then follows from an envelope 
property   that also leads to the variational representation of the limiting height profile.
In this section we treat the wedge case, and the next section puts it all together.  

Recall the notation and conventions introduced in the previous section.  In particular, 
$c(x)$ is  a positive, lower semicontinuous speed function with only   finitely many discontinuities in any compact set. Define a lattice analogue of the wedge $\cW$ by 
\be\cL=\{(i,j) \in \bZ^2: j\geq1, i\geq -j+1\}\label{defLL}\ee
 with boundary $\partial\cL=\{(i,0):i\geq 0\}\cup\{(i,-i): i<0\}$.

For each $n\in \bN$   construct a last-passage growth model on $\cL$ that represents the TASEP height function in the wedge. Let $\{ \tau^n_{i,j}: (i,j) \in \cL \}_{n\in \bN}$ denote a sequence of independent collections of i.i.d.\ exponential rate $1$ random variables. We need an extra index $\ell$ to denote the shifting. Define weights
\be
\omega_{i,j}^{n,\ell} = c\left(\dfrac{i - \ell}{n}\right)^{-1}, \quad (i,j)\in \mathcal{L}. 
\label{weights}
\ee 
For $\ell \in \bZ$ and  $n \in \bN$ assign to   site $(i,j) \in \cL$ the random variable $\omega_{i,j}^{n,\ell}\tau^n_{i,j}.$ 
Given lattice points $(a,b), (u,v) \in \cL$,   $\Pi((a,b),(u,v))$ is the set of lattice paths $\pi = \{(a,b) = (i_0,j_0),(i_1,j_1),...,(i_p,j_p) = (u,v)\}$ whose admissible steps satisfy 
\be
(i_l,j_l) - (i_{l-1},j_{l-1}) \in \{ (1,0), (0,1), (-1,1)\}.
\label{wedgepaths}
\ee
In the case that $(a,b)=(0,1)$ we simply denote this set by $\Pi(u,v)$. 
For $(u,v)\in \cL$, $\ell \in \bR$ and $n\in \bN$ denote the \textsl{wedge last passage time}
\be
T^{n,\ell}(u,v) = \max_{\pi \in \Pi(u,v)} \sum_{(i,j)\in \pi}\omega^{n,\ell}_{i,j}\tau^n_{i,j}
\label{lastpassagetime}
\ee with boundary conditions
\be
T^{n,\ell}(u,v)=0 \quad \text{for}\quad (u,v)\in \partial\cL.
\label{boundaryoflastpassagetime}
\ee

Admissible 
steps \eqref{wedgepaths}  come from the properties of the TASEP height function.  
Notice that $(0,1)$ is in fact  never used in a maximizing path.   

To describe macroscopic last passage times define,  for $(x,y)\in \cW$ and $q \in \bR$, 
\be
\Gamma^q(x,y)=\sup_{\textbf{x}(\cdot)\in \cH(x,y)}\Big\{ \int_{0}^{1}\frac{\gamma(\textbf{x}'(s))}{c(x_1(s)-q)}\,ds \Big\}.
\label{gammaq}
\ee
 
 \begin{theorem}
For all $q\in \bR$ and $(x,y)$ in the interior of $\cW$
\be
\dlim n^{-1}T^{n,\fl{nq}}(\fl{nx},\fl{ny})=\Gamma^q(x,y) \quad \text{a.s.}
\label{Tlimeq} \ee\label{Tlimit}
\end{theorem} 

\begin{remark} In a constant 
 rate $c$ environment the wedge last passage limit is
\be
\lim_{n\rightarrow\infty}\frac1{n} {T^{n}(\fl{nx},\fl{ny})} = c^{-1}\gamma(x,y)  = c^{-1}\left(\sqrt{x+y}+\sqrt{y}\right)^2.
\label{chomogeneous}
\ee 
The limit  $\gamma(x,y)$
  is   concave, but this is  not true 
 in general for   $\Gamma^0(x,y)$.
In some special cases concavity still holds,
such as   if the function $c(x)$ is nonincreasing if $x<0$ and nondecreasing if $x>0$. 
\end{remark} 

To prove Theorem \ref{Tlimit} we approximate $c(x)$ with step   functions. Let $-\infty = a_1 <a_2< ... < a_{L-1} <a_L = +\infty,$ and consider the lower semicontinuous step function 
\be 
c(x) = \sum_{m=1}^{L-1}r_m \mathbf{1}_{(a_m, a_{m+1})}(x) + \sum_{m=2}^{L-1}\min\{r_{m-1}, r_{m}\}\mathbf{1}_{\{ a_m \}}(x).
\label{simpleratefunction}
\ee    

\begin{proposition}
Let $c(x)$ be given by \eqref{simpleratefunction}. Then limit \eqref{Tlimeq} holds. 
 \label{simpleTlimit}
\end{proposition}

On the way to Proposition \ref{simpleTlimit} we state preliminary lemmas that will be used for
pieces of paths.  We write $c_i$ for the rate values instead of $r_i$ 
to be consistent with the notation in Theorem \ref{two-phaselastpassagetime}. 

\begin{lemma}
Assume that there is a unique discontinuity $a_2 = 0$ for the speed function $c(x)$ in \eqref{simpleratefunction}. Then for $y>0$
\[
\dlim n^{-1}T^{n,0}(0,\fl{ny}) = \dfrac{4y}{\min\{c_1, c_2\}} = \int_{0}^{1}\frac{\gamma(0,y)}{c(0)}\,ds\quad \text{a.s.} 
\]
\label{Verticalpassage}
\end{lemma}

\begin{proof} The upper bound in the limit is immediate from domination with  
constant rates $c(0)$. 

 For the lower bound we spell out the details for the case
 $c_1 \geq c_2$.
  Let $\varepsilon >0$.  To bound $T^{n,0}(0,\fl{ny})$ from below 
  force  the  path  to  go through  points  
$(0,1)$, $ \{(\fl{ny\e}, (k-1)\fl{ny\e}):  k=1,\dotsc, \fl{\e^{-1}}\}$ and  $(0,\fl{ny})$. 
For  $1\le k< \fl{\e^{-1}}$  let $T^n(R^n_k)$ be the last passage time from $(\fl{ny\e}, (k-1)\fl{ny\e})$  to 
  $(\fl{ny\e}, k\fl{ny\e})$.   $R^n_k$ refers to  the parallelogram that contains all the
  admissible paths   from $(\fl{ny\e}, (k-1)\fl{ny\e})$  to 
  $(\fl{ny\e}, k\fl{ny\e})$.  Each $R^n_k$ lies to the right of $x=0$ and therefore 
   in the $c_2$-rate area.  (See Fig.~\ref{lowrectangles}.) 
  
\begin{figure}[ht]
\begin{center}
\begin{picture}(220,220)(-100,0)
\put(20,20){\vector(1,0){200}} %axis in the direction of (1,0), length 160 pixels 
\put(20,20){\vector(-1,1){200}} %axis in the direction of (0,1) 
\put(20,20){\vector(0,1){200}} %axis in the direction of (0,1) 

%Square with diagonal
\put(-160,200){\line(1,0){180}}
\put(200,20){\line(-1,1){180}}
\put(20,20){\line(0,1){180}}

%Rectangles
\put(70,20){\line(-1,1){50}}
\put(120,20){\line(-1,1){100}}
\put(20,70){\line(1,0){60}}
\put(20,120){\line(1,0){10}}
\put(60,120){\line(1,0){40}}
\put(110,70){\line(1,0){10}}
\put(120,70){\line(-1,1){100}}

%Path 
{\color{nicosred}
\linethickness{1.1 pt}
\put(20,20){\line(1,0){50}}
\put(70,20){\line(-1,1){10}}
\put(60,30){\line(1,0){9}}
\put(69,30){\line(-1,1){10}}
\put(59,40){\line(1,0){5}}
\put(64,40){\line(-1,1){10}}
\put(54,50){\line(1,0){12}}
\put(66,50){\line(-1,1){10}}
\put(56,60){\line(1,0){14}}
\put(70,60){\line(-1,1){10}}
\put(60,70){\line(1,0){14}}
\put(74,70){\line(-1,1){14}}
\put(60,84){\line(1,0){4}}
\put(64,84){\line(-1,1){10}}

\put(70,115){\line(1,0){5}}
\put(75,115){\line(-1,1){12}}
\put(63,127){\line(1,0){6}}
\put(69,127){\line(-1,1){12}}
\put(57,139){\line(1,0){8}}
\put(65,139){\line(-1,1){10}}
\put(55,149){\line(1,0){10}}
\put(65,149){\line(-1,1){10}}
\put(55,159){\line(1,0){4}}
\put(59,159){\line(-1,1){10}}
\put(49,169){\line(1,0){2}}
\put(51,169){\line(-1,1){31}}

}

%Thickens Path 
{\color{nicosred}
%\linethickness{1.5 pt}
%\put(17,20){\line(1,0){50}}
\put(67,20){\line(-1,1){10}}
%\put(57,30){\line(1,0){9}}
\put(65,30){\line(-1,1){10}}
%\put(59,40){\line(1,0){5}}
\put(61,40){\line(-1,1){10}}
%\put(54,50){\line(1,0){12}}
\put(63,50){\line(-1,1){10}}
%\put(56,60){\line(1,0){14}}
\put(67,60){\line(-1,1){10}}
%\put(60,70){\line(1,0){14}}
\put(71,70){\line(-1,1){14}}
%\put(60,84){\line(1,0){4}}
\put(61,84){\line(-1,1){10}}

%\put(70,115){\line(1,0){5}}
\put(72,115){\line(-1,1){12}}
%\put(63,127){\line(1,0){6}}
\put(66,127){\line(-1,1){12}}
%\put(57,139){\line(1,0){8}}
\put(62,139){\line(-1,1){10}}
%\put(55,149){\line(1,0){10}}
\put(62,149){\line(-1,1){10}}
%\put(55,159){\line(1,0){4}}
\put(56,159){\line(-1,1){10}}
%\put(49,169){\line(1,0){2}}
\put(48,169){\line(-1,1){31}}

}

%Dots
\put(57,100){\vdots}
%Strips
{\color{light-gray}
\put(7,30){\line(0,1){170}}
\put(-3,40){\line(0,1){160}}
\put(-13,50){\line(0,1){150}}
\put(-23,60){\line(0,1){140}}
\put(-33,70){\line(0,1){130}}
\put(-43,80){\line(0,1){120}}
\put(-53,90){\line(0,1){110}}
\put(-63,100){\line(0,1){100}}
\put(-73,110){\line(0,1){90}}
\put(-83,120){\line(0,1){80}}
\put(-93,130){\line(0,1){70}}
\put(-103,140){\line(0,1){60}}
\put(-113,150){\line(0,1){50}}
\put(-123,160){\line(0,1){40}}
\put(-133,170){\line(0,1){30}}
\put(-143,180){\line(0,1){20}}
\put(-153,190){\line(0,1){10}}
}

%Coordinate system x-axis
\put(10,10){\small 0}
\put(60,10){$\fl{\e y n}$}

{\color{nicosred}\put(80,30){\small{$R^n_1$}}}

\put(110,10){$2\fl{\e y n}$}
\put(160,10){$\cdots$}

{\color{nicosred}\put(90,130){\small{$R^n_{\fl{\e^{-1}}}$}}
\put(90,135){\vector(-1,0){15}}
}

\put(190,10){$\fl{ y n}$} 
%Coordinate system y-axis
\put(-65,62){$\fl{\e y n}$}
\put(-115,112){$2\fl{\e y n}$}
%\put(-175,155){$\vdots$}
\put(-188,194){$\fl{ y n}$} 

%Dashed lines
\multiput(20,70)(-5,0){11}{\line(-1,0){2}}
\multiput(20,120)(-5,0){21}{\line(-1,0){2}}
\multiput(17,170)(-5,0){30}{\line(-1,0){2}}
\end{picture}
\end{center}
\caption{A possible microscopic path   forced to go through opposite corners of the parallelograms $R^n_k$. The striped area left of   $x=0$ is the $c_1$-rate region. }
\label{lowrectangles}
\end{figure}
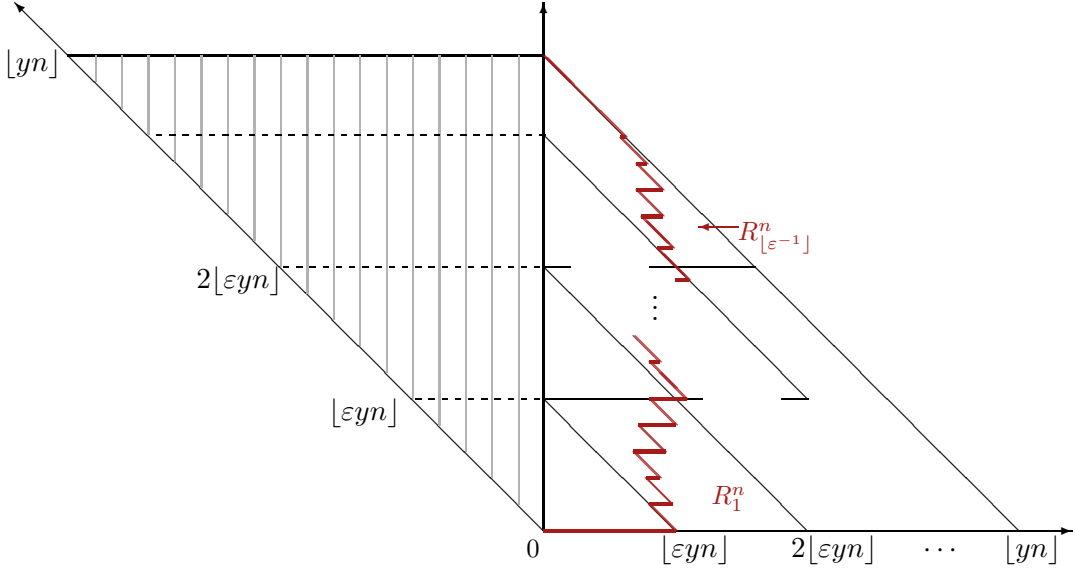

Let $0< \delta < \e {c_2}^{-1}\gamma(0,y)$. 
A large deviation estimate (Theorem 4.1 in \cite{sepp-large-deviations}) 
gives  a constant $C = C(c_2,y,\e,\delta)$ such that 
\be
\bP\left\{ T^n_{c_2}(R^n_k)\leq n(\e {c_2}^{-1}\gamma(0,y) - \delta) \right\}\leq e^{-Cn^2}.
\label{largedeviations}
\ee
 By a   Borel-Cantelli argument, for large $n$, 
\begin{align*}
T^{n,0}(0,\fl{ny}) &\ge \sum_{k=1}^{\fl{\e^{-1}}-1}T^n(R^n_k) \ge n(\fl{\e^{-1}}-1)
(\e {c_2}^{-1}\gamma(0,y) - \delta).  
\end{align*}
This suffices for the conclusion.  
 \end{proof}
\begin{remark} 
This lemma shows why it is convenient to use a lower semi-continuous speed function. A path that starts and ends at the same discontinuity stays mostly in the low rate region to maximize its weight. This translates macroscopically to the formula for the limiting time constant obtained in the lemma, involving only the value of $c$ at the discontinuity. If the speed function is not lower semi-continuous, we can state the same result using left and right limits.
\end{remark}

\begin{lemma}
Let  $ a = 0 < b< +\infty$
be discontinuities  for the step speed function $c(x)$ and $c(x)=r$ for $a<x<b$.  
Take  $z\in[0,b]$.  Let 
$\wt T^n(\fl{nz}, \fl{ny})$ be the wedge last passage time from $(0,1)$ to $(\fl{nz},\fl{ny})$ subject to the constraint that the   path has to stay in the $r$-rate region $(a,b)\times(0,+\infty)$,
except possibly for the initial and final steps. Then  
\be
\dlim n^{-1}\wt T^n(\fl{nz},\fl{ny}) = r^{-1}{\gamma(z,y)} \quad \text{a.s.} 
\ee
Same statement holds if $b\le z\le 0$.  
\label{Conditionalpassage}
\end{lemma}

\begin{proof}
The upper bound  $\varlimsup n^{-1}\wt T^n(\fl{nz},\fl{ny}) 
\leq {r}^{-1}{\gamma(z,y)}$  is immediate
by  putting constant rates $r$ everywhere and dropping the restrictions
  on the path. For the lower bound
adapt   the steps of the proof in Lemma \ref{Verticalpassage}. \end{proof}

Lemma \ref{Conditionalpassage} is a place where we cannot allow accumulation of discontinuities for the speed function.   

\begin{figure}[ht]
\begin{center}
\begin{picture}(150,90)(20,0)
\put(20,20){\vector(1,0){175}} %axis in the direction of (1,0), length 160 pixels 
\put(80,20){\vector(0,1){80}} %axis in the direction of (0,1) 
%Parallelograms
\put(80,20){\line(-1,1){60}} 
\put(185,20){\line(-1,1){60}}
\put(20,80){\line(1,0){100}}

%Discontinuity
{\color{darkblue}
\linethickness{1.5 pt}
\put(125,20){\line(0,1){60}}
\put(80,20){\line(0,1){60}}
}

%Other Regions
%right
{\color{light-gray}
\put(125,20){\line(0,1){55}}
\put(130,20){\line(0,1){50}}
\put(135,20){\line(0,1){45}}
\put(140,20){\line(0,1){40}}
\put(145,20){\line(0,1){35}}
\put(150,20){\line(0,1){30}}
\put(155,20){\line(0,1){25}}
\put(160,20){\line(0,1){20}}
\put(165,20){\line(0,1){15}}
\put(170,20){\line(0,1){10}}
\put(175,20){\line(0,1){5}}
}

%Other Regions
%left
{\color{light-gray}
\put(70,25){\line(0,1){55}}
\put(65,30){\line(0,1){50}}
\put(60,35){\line(0,1){45}}
\put(55,40){\line(0,1){40}}
\put(50,45){\line(0,1){35}}
\put(45,50){\line(0,1){30}}
\put(40,55){\line(0,1){25}}
\put(35,60){\line(0,1){20}}
\put(30,65){\line(0,1){15}}
\put(25,70){\line(0,1){10}}
\put(20,75){\line(0,1){5}}
}
%Path 
{\color{nicosred}
\put(72,20){\line(1,0){25}}
\put(97,20){\line(-1,1){10}}
\put(87,30){\line(1,0){8}}
\put(95,30){\line(-1,1){10}}
\put(85,40){\line(1,0){14}}
\put(99,40){\line(-1,1){10}}
\put(89,50){\line(1,0){16}}
\put(105,50){\line(-1,1){10}}
\put(95,60){\line(1,0){8}}
\put(103,60){\line(-1,1){10}}
\put(93,70){\line(1,0){17}}
\put(110,70){\line(-1,1){10}}
\put(100,80){\line(1,0){15}}
}
%Coordinate system x-axis
\put(60,10){\small 0}

\end{picture}
\end{center}
\caption{A possible microscopic path described in Lemma \ref{Conditionalpassage}. The path has to stay in the unshaded region. }
\label{paths}\end{figure}
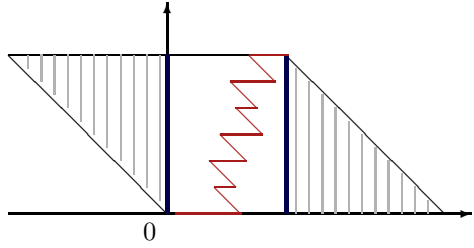

  Before proceeding to the proof of Proposition \ref{simpleTlimit} we make
  a simple but important  observation  about the macroscopic paths 
  $\mathbf{x}(s)= (x^1(s), x^2(s))$, $s\in[0,1]$, in $\cH(x,y)$   
  for the case where  $c(x)$ is a step function \eqref{simpleratefunction}.   
    
\begin{lemma}
There exists a constant $C=C(x,y,c(\cdot-q))$ such that
the supremum in \eqref{gammaq} comes from 
 paths in  $\cH(x,y)$ that consist of 
at most $C$  line segments.  Apart from the first and last
segment, these segments can be of two types:  segments that 
go from one discontinuity of $c(\cdot-q)$ to a neighboring 
discontinuity, and vertical segments along a discontinuity. 
\label{segments}
\end{lemma}

 \begin{proof}   
Path $\mathbf x$  is  a union of subpaths $\{\bx_j\}$  along which  
$c(x_j^1(s)-q)$ 
is constant, except possibly at the endpoints.
  % (An extreme case  or is a vertical line segment that lies on a discontinuity of $c(x)$. 
  Given such a subpath $(\mathbf x_j(s): {t_j}\le s\le {t_{j+1}})$,  concavity of $\gamma$ and Jensen's inequality imply that the   line segment  
  $\phi_j$ that connects $\bx_j(t_j)$ to  $\bx_j(t_{j+1})$ dominates:  
\be
\int_{t_j}^{t_{j+1}} \frac{\gamma(\bx_j '(s))}{c(x_j^1(s)-q)}\,ds \leq \int_{t_j}^{t_{j+1}} \frac{\gamma(\phi_j '(s))}{c(\phi_j^1(s)-q)}\,ds.
\notag 
\ee
Consequently we can restrict to paths that are unions of line segments.  

To bound the number of line segments, observe first that 
the number of segments that go from one discontinuity to a neighboring
discontinuity is bounded.  The reason is that the restriction
$\bx'(s)\in\cW$ forces  such a segment  to
increase at least one of the coordinates by the distance between
the discontinuities.  

Additionally there can be subpaths that  touch
 the same discontinuity more than once 
without touching a different discontinuity. Lower semi-continuity
of $c(\cdot)$  and Jensen's inequality
show again that the vertical line segment  that stays on the discontinuity 
dominates such a subpath.  
Consequently there can be at most one (vertical) line segment 
between two line segments that connect distinct discontinuities.  
\end{proof} 
 
Next a lemma about the continuity of $\Gamma^q$. We write 
$\Gamma^q((a,b),(x,y))$ for the value in \eqref{gammaq} 
when the paths go from $(a,b)$ to $(x,y)\in (a,b)+\cW$. 

\begin{lemma}
Fix $z, w >0$. Then there exists 
%$\delta_* = \delta_*(z,w, c(\cdot-q))>0 $ and 
a constant $C= C(z,w, c(\cdot-q))<\infty$ such  that for all 
$0< \delta \leq 1$ and $0 \leq a \leq z$
\be
\Gamma^q((a,0),(z,\delta)) - \Gamma^q((a,0),(z,0)) \leq C\sqrt{\delta},
\label{deltaparineq1}\ee
and for $0 \leq b \leq w$
\be
\Gamma^q((-b,b),(-w,w+\delta)) - \Gamma^q((-b,b),(-w,w)) \leq C\sqrt{\delta}.
\label{deltaparineq2}\ee
\label{deltapar}
\end{lemma}

\begin{proof}
Pick $\delta\in(0,1]$ and consider the point
 $(z,\delta)$ in $\cW$.  For any $\mathbf{x}= (x^1(s),x^2(s)) \in \cH(z,\delta)$ set 
\be 
I(\bx,q) =\int_{0}^{1}\dfrac{\gamma(\bx'(s))}{c(x^1(s)-q)}\,ds.  
\ee 
Let $\e>0 $ and assume that $\phi=(\phi^1,\phi^2)\in \cH(z,\delta)$ is a path such that $\Gamma^q(z,\delta) - I(\phi,q) < \e.$
 Lemma \ref{segments} implies that we can decompose $\phi$ into 
 disjoint linear segments $\phi_j$ so that $\phi= \sum_{j=1}^M \phi_j$
and $\phi_j:[s_{j-1},s_{j}]\to\cW$. Here $\sum_{j}\phi_j$ means path concatenation.

We can find segments $\phi_{j(k)}$, $1\le k\le N$,  such that 
\[\phi^1_{j(k)}(s_{j(k)-1})< \phi^1_{j(k)}(s_{j(k)}), \quad  
\phi^1_{j(k)}(s_{j(k)})= \phi^1_{j(k+1)}(s_{j(k+1)-1}), \]
$\phi^1_{j(1)}(s_{j(1)-1})=0$, and $\phi^1_{j(N)}(s_{j(N)})=z$. 
 In other words, the projections of the  segments  $\phi_{j(k)}$ cover the interval
 $[0, z]$ without overlap and without backtracking.  

We bound the contribution of the remaining path segments to $I(\phi, q)$. 
Let $J$ $=$ \break $\bigcup_{k=1}^{N-1} [s_{j(k)}, s_{j(k+1)-1}]$ be the 
leftover portion of the time interval $[0,1]$.  
The subpath $\phi(s)$, $s\in [s_{j(k)}, s_{j(k+1)-1}]$,
  (possibly) eliminated from between 
$\phi_{j(k)}$ and $\phi_{j(k+1)}$  satisfies 
$\phi^1(s_{j(k)})=\phi^1(s_{j(k+1)-1})$.   
Note that $\gamma(a,b)\le 2a+4b$ for $(a,b)\in\cW$ and 
$\int_0^1 (\phi^2)'(s)\,ds = \delta$. 
  We
 can bound as follows: 
  \begin{align}
  \int_{J}\frac{\gamma((\phi^1)'(s),(\phi^2)'(s))}{c(\phi^1(s)-q)}\,ds 
       &\leq C\int_{J}\gamma((\phi^1)'(s),(\phi^2)'(s))\,ds \notag\\
       &\leq C\int_{J} \big(2(\phi^1)'(s)+4(\phi^2)'(s)\big)\,ds \notag\\
       &\leq C\int_{J} 2(\phi^1)'(s)\,ds +  C\int_{0}^{1}4(\phi^2)'(s)\,ds \notag\\
       &= 0 + 4C\delta.       
  \end{align}
 
%For $1\leq j \leq N$ and consider the segment $\phi_{j(k)}$. 
Set $t_k = s_{j(k)-1} < u_k = s_{j(k)}$. 
%Let $\phi_{j(k)}(t_j) = (a_j, b_j)$, 
%$\phi_{j(k)}(t_{j+1}) = (a_{j+1},b_{j+1})$, 
%where $a_j < a_{j+1}$ are discontinuities of $c(\cdot - q).$ 
Define a horizontal path $w$ from $(0,0)$ to $(z,0)$ 
with segments 
\be
w_{k}(s)= \big(\phi_{j(k)}^{1}(s), 0\big), \quad \text{for } 
\  t_k \leq s \leq u_k, 
\ee
and constant on the complementary time set $J$.

To get the lemma, we estimate \normalsize{
\begin{align}
\Gamma^q(z&,\delta) -\e \leq I(\phi,q) = \int_J \frac{\gamma(\phi'(s))}{c(\phi^1(s)-q)}\,ds +  \int_{[0,1]\setminus J} \frac{\gamma(\phi'(s))}{c(\phi^1(s)-q)}\,ds\notag \\ 
                   &\leq C\delta + \sum_{k = 1}^{N} \Big(I(\phi_{j(k)},q) - I(w_{k}, q)\Big)+ \Gamma^q(z,0) \notag \\
                   &\leq C\delta + C'\sum_{k=1}^{N}\int_{t_{k}}^{u_{k}} \big(\gamma(\phi'_{j(k)}(s)) - \gamma(w'_{k}(s))\big)\,ds+  \Gamma^q(z,0)\notag \\
                   &\leq C\delta + C'\sum_{k=1}^{N}\bigg(\int_{t_{k}}^{u_k} (\phi^2)'_{j(k)}(s)\,ds +\notag \\ &\phantom{xxxxxxxxxxxxxxx}+2\int_{t_{k}}^{u_k}\sqrt{(\phi^2)'_{j(k)}(s)}\sqrt{(\phi^1)'_{j(k)}(s)+(\phi^2)'_{j(k)}(s)} \,ds\bigg)+\Gamma^q(z,0)\notag \\ 
                    &\leq C\delta + C'\sum_{k=1}^{N}\bigg(\int_{t_{k}}^{u_k} (\phi^2)'_{j(k)}(s)\,ds\bigg)^{\frac{1}{2}}\bigg(\int_{t_{k}}^{u_k}\big((\phi^1)'_{j(k)}(s)+(\phi^2)'_{j(k)}(s)\big) \,ds\bigg)^{\frac{1}{2}}+  \Gamma^q(z,0)\notag \\ 
                    &\leq C\delta + C'\bigg(\sum_{k=1}^{N}\int_{t_{k}}^{u_k} (\phi^2)'_{j(k)}(s)\,ds\bigg)^\frac 12\times \notag \\
                    &\phantom{xxxxxxxxxxxxx}\times\bigg(\sum_{k=1}^{N}\int_{t_{k}}^{u_k}\big((\phi^1)'_{j(k)}(s)+(\phi^2)'_{j(k)}(s)\big) \,ds\bigg)^{\frac{1}{2}}+  \Gamma^q(z,0)\notag \\ 
                    &\leq C\delta + C'\sqrt{\delta}\sum_{k=1}^{N}\int_{t_{k}}^{u_k}\big((\phi^1)'_{j(k)}(s)+(\phi^2)'_{j(k)}(s) \big)\,ds+\Gamma^q(z,0)\notag \\ 
                    &\leq C\delta + C'\sqrt{\delta}\sqrt{z+\delta}+  \Gamma^q(z,0)\notag \\
                    &\leq C\delta + C'\sqrt{\delta}\sqrt{z}+ \Gamma^q(z,0).\notag  
%\label{horizontalestimate}        
\end{align}}
 The first
 inequality \eqref{deltaparineq1} follows
 for $a = 0$ by letting $\e$ go to 0.  It also follows
 for all $a\in[0, z]$ by shifting the origin to $a$ which 
replaces $z$ with $z-a$.

For the second inequality \eqref{deltaparineq2} the arguments are 
analogous, so we omit them.
\end{proof}

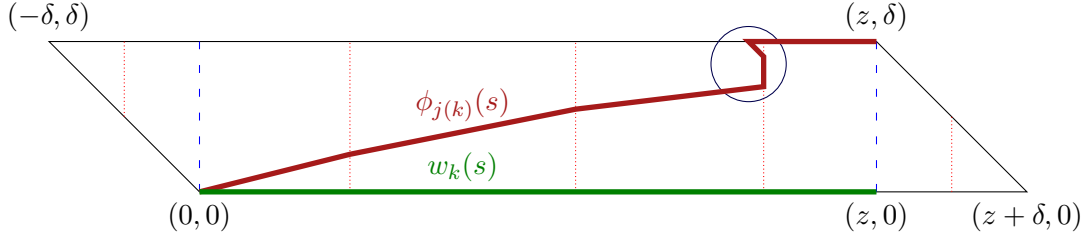
\begin{figure}
\begin{tikzpicture}
\draw (0,0) -- (11,0)--(9,2)--(-2,2)--(0,0);
\draw (0,0) node[below]{$(0,0)$};
\draw (-2,2) node[above]{$(-\delta,\delta)$};
\draw (9,2) node[above]{$(z,\delta)$};
\draw (9,0) node[below]{$(z,0)$};
\draw (11,0) node[below]{$(z+\delta, 0)$};

\draw[densely dotted, red] (-1,1) -- (-1,2) ;
\draw[densely dotted, red]  (2,0)--(2,2);
\draw[densely dotted, red]  (5,0)--(5,2);
\draw[densely dotted, red]  (7.5,0)--(7.5,2);
\draw[densely dotted, red]  (10,0)--(10,1);

\draw[color=darkblue] (7.3, 1.7) circle (.5);

\draw[loosely dashed, blue](0,0)--(0,2);
\draw[loosely dashed, blue](9,0)--(9,2);

\draw[nicosred, line width=2pt] (0,0)--(2,.5)-- (5,1.1)--(7.5,1.4)--(7.5, 1.8)--(7.3,2)--(9,2);
\draw[nicosred] (3.5,.8) node[above]{$\phi_{j(k)}(s)$};

\draw[darkgreen, line width=2pt] (0,0)--(9,0);
\draw[darkgreen](3.5,0) node[above]{$w_k(s)$};

\end{tikzpicture}
\caption{ A possible macroscopic path from $(0,0)$ to $(z,\delta)$. 
 The dotted vertical lines  are discontinuity columns of $c(\cdot-q).$ The error from eliminating the segments outside the two vertical dashed lines and from eliminating  pathologies (like the circled part) is of order $\delta$ and a comparison with the horizontal path leads to an error of order $\sqrt{\delta}$.}
\end{figure}

\begin{corollary}
Fix $(x,y) \in \cW$.  Then there exists  $C=C(x, y, c(\cdot - q))<\infty$
such that for all
 $0<\delta \le 1 $ 
\be
\Gamma^q(x,y+\delta) - \Gamma^q(x,y) < C\sqrt{\delta}.
\ee
\label{deltapar2}
\end{corollary}

\begin{proof}
Let $A((a,b),(x,y))$ be the parallelogram with sides parallel to
 the boundaries of the wedge, north-east corner the point $(x,y)$ and 
south-west corner at $(a,b)$. If $(a,b) = (0,0)$ we simply write $A(x,y)$.

Let $\e > 0.$ Let $\phi^{\e}$ a path such that $\Gamma^q(x,y+\delta) - I(\phi^{\e}, q)< \e$. Let $u$ be the point where $\phi^{\e}$ 
first intersects the north or the east  boundary of $A(x,y)$.
 Without loss of generality assume it is the north boundary
and so  $u = (a,y)$ for some $a\in[-y,x]$.  Then,
\begin{align}
\Gamma^q(x,y+\delta)-\e &\leq I(\phi^{\e}, q) \notag \\
                        &\leq \Gamma^q(a,y)+\Gamma^q((a,y),(x,y+\delta)) \notag \\
                        &=\Gamma^q(a,y)+\Gamma^q((a,y),(x,y)) + \Gamma^q((a,y),(x,y+\delta))-\Gamma^q((a,y),(x,y))\notag \\
                        &\leq \Gamma^q(x,y)+ \Gamma^q((a,y),(x,y+\delta))-\Gamma^q((a,y),(x,y)). 
\end{align}
The last inequality gives
\be
\Gamma^q(x,y+\delta) - \Gamma^q(x,y)\leq \Gamma^q((a,y),(x,y+\delta))-\Gamma^q((a,y),(x,y)) +\e \leq C\sqrt{\delta} +\e
\ee
by Lemma \ref{deltapar}. Let $\e$ decrease to $0$ to prove the Corollary.
\end{proof}

\begin{proof}[Proof of Proposition \ref{simpleTlimit}]
Fix $(x,y)$ in the interior of $\cW$.  
For   $\mathbf{x}= (x^1(s),x^2(s)) \in \cH(x,y)$ set 
\be 
I(\bx,q) =\int_{0}^{1}\dfrac{\gamma(\bx'(s))}{c(x^1(s)-q)}\,ds.  
\ee 

We prove first 
\be
\varliminf_{n\rightarrow \infty} n^{-1}T^{n,\fl{nq}}(\fl{nx},\fl{ny})\geq \Gamma^q(x,y)
\equiv \sup_{\textbf{x}(\cdot)\in \cH(x,y)} I(\bx,q).    
\ee
It suffices to consider  macroscopic paths   of the type  
\be
\bx(s) = \sum_{j=1}^{H}\bx_j(s)\mathbf{1}_{[s_j,s_{j+1})}(s)
\label{type22}  \ee 
where $H\in \bN$,   $\bx_j$ is the straight line segment 
from  $\bx(s_j)$ to  $\bx(s_{j+1})$,    $c(x_1(s) - q) = r_{m_j}$ is constant for $s\in( s_j, s_{j+1})$, 
and by continuity   $\bx_j(s_{j+1}) = \bx_{j+1}(s_{j+1})$.

 Let $\pi^n$ be the  microscopic path through points 
$(0,1)$,   $ \{\fl{n\mathbf x_j(s_j)}: 1\le j\le K\}$ and  $(\fl{nx},\fl{ny})$
constructed so that its segments $\pi^n_j$  satisfy these requirements: 

(i)   $\pi^n_j$   lies inside the region where $\omega^{n,\fl{nq}}_{i,k} = r_{m_j}^{-1}$ 
is constant, except possibly  for the initial and final step;   

(ii) $\pi^n_j$ maximizes passage time between its endpoints $\fl{n\textbf{x}_j(s_j)}$
and  $ \fl{n\bx_{j+1}(s_{j+1})}$ subject to the above requirement.

Let 
\be T^{n,\fl{nq}}_j = \max_{\pi^n_j} \sum_{(i,k) \in \pi^n_j}\omega^{n,\fl{nq}}_{i,k}\tau_{i,k}
\label{Tpieces}
\ee 
denote the  last-passage time of a segment subject to these constraints.  
Observe that the proofs of  Lemmas \ref{Verticalpassage} and  \ref{Conditionalpassage} do not depend on the shift parameter $q$, therefore
\begin{align*}
\dlim n^{-1}T^{n,\fl{nq}}_j = \dfrac{\gamma( \bx_j(s_j) - \bx_{j+1}(s_{j}))}{r_{m_j}}  
%=\dfrac{\gamma(\bx_j(s_j) - \bx_{j}(s_{j+1}))}{c(x^1_j(s_j^+) - q)} \notag \\
                                   = \int_{s_j}^{s_{j+1}}\dfrac{\gamma( \bx'_j(s))}{c(x^1_j(s) - q)}\,ds.
\end{align*}
Adding up the segments gives the lower bound: 
\begin{align*}
\varliminf_{n\rightarrow \infty} n^{-1}T^{n,\fl{nq}}(\fl{nx},\fl{ny})
&\ge  \varliminf_{n\rightarrow \infty} \sum_{j} n^{-1} T^{n,\fl{nq}}_j \\
 &=\sum_{j} \int_{s_j}^{s_{j+1}}\dfrac{\gamma(\bx'_j(s))}{c(x^1(s) -q)}\,ds  
                                                       =I(\bx,q).   
\end{align*} 

Now for the complementary upper bound 
\be
\varlimsup_{n\rightarrow \infty} n^{-1}T^{n,\fl{nq}}(\fl{nx},\fl{ny})\leq \Gamma^q(x,y).
\label{ub7}\ee
Each microscopic path  to $(\fl{nx},\fl{ny})$ 
 is contained in $nA$ for a fixed parallelogram 
$A \subseteq \mathcal{W}$ with sides parallel to the wedge boundaries. 
Pick $\e>0$.
 Let $r_*>0$ be a lower bound on all the rate values that appear in the set $A$. 
Find $\delta>0$ such that   $\abs{\gamma(v) - \gamma(w)}  < \e r_*$ 
 for all $v,w \in A$ with $\abs{v-w} < \delta$ and $\delta \le 1$ 
 so that Corollary \ref{deltapar2} is valid.

Consider an arbitrary microscopic path from $(0,1)$ to 
$(\fl{nx},\fl{ny})$.   Given the speed function and $q$,
there is a fixed upper bound $Q=Q(x,y)$ 
on the number of segments of the path that start at one
discontinuity column $(\fl{na_i} + \fl{nq})\times\bN$ and 
end at a neighboring discontinuity column 
$(\fl{na_{i\pm 1}} + \fl{nq})\times\bN$. The reason
is that  there is an
order $n$ lower bound on the number of lattice steps it takes 
to travel between distinct discontinuities  
in $nA$.

 Fix $K\in\bN$ and partition the interval $[0,y]$ evenly by $b_j=jy/K$, $0\le j\le K$,  
so that $y/K < \delta/Q.$ Make the partition finer by adding 
the $y-$coordinates of the intersection points of discontinuity lines
$\{a_i+q\}\times\bR_+$  with the boundary of $A$. 

Let $\pi^n$ be the maximizing microscopic path.  
We decompose   $\pi^n$ into  path 
segments  $\{\pi^n_{j}: 0\le j<M_n\}$
by looking at visits to   discontinuity columns $(\fl{na_i} + \fl{nq})\times\bN$, 
both repeated visits to the same
column and visits to a column different from the previous one.  
Let  
$\{0=b_{k_0}\leq b_{k_1} \leq b_{k_2} \leq ... \leq b_{k_{M_n-1}} \leq b_{k_{M_n}}= y\}$ 
be a sequence of partition points  and   
 $\{ 0=x_0,\, x_1=a_{m_1}+q,\, x_2=a_{m_2}+q, \dotsc, x_{M_n}=x\}$ 
a sequence where $x_{j}$ for $0<j<M_n$ are discontinuity points of the 
shifted speed function $c(\cdot\,-q)$. 
We can create the path segments and these sequences  
with the property that 
%for $0<j<M_n-1$ the 
segment $\pi^n_{j}$ starts at  
  $(\fl{nx_{j}}, l)$ with $l$ in the range $\fl{nb_{k_j}} \leq l \leq \fl{nb_{k_{j}+1}}$ and ends at $(\fl{nx_{j+1}}, l')$ with $\fl{nb_{k_{j+1}}} \leq l' \leq \fl{nb_{k_{j+1}+1}}.$
%
%  Additionally there are initial and final segments $\pi^n_0$ and 
%$\pi^n_{M_n-1}$  that possibly have only one endpoint at a discontinuity column. 
In an extreme
case the entire path $\pi^n$ can be a single segment that does  not
touch discontinuity columns.  

In order to have a  fixed  upper bound on the total number 
$M_n$ of   segments, uniformly in $n$, we insist that
for $0<j<M_n-1$ the labels satisfy:

(i) For odd $j$, $\pi^n_j$ starts and ends at the same discontinuity column $(\fl{nx_{{j}}}, \,\cdot\,)$.    The rate relevant for segment $\pi^n_j$
is  $r_{\ell_j}=c(a_{m_j})$. 

(ii) For  even $j$, $\pi^n_j$ starts and ends at different neighboring
discontinuity columns, and  except for the initial and final points,
  does not touch any discontinuity column
 and visits only points that are in a 
region of constant rate $r_{\ell_j}$.  

The above conditions may create empty segments.  This is not harmful.
Replace $Q$ with $2Q+2$ to continue having the uniform upper bound
 $M_n\le Q$.  

Let $T(\pi^n_j)$  be the total weight of segment $\pi^n_j$.
Let $\tilde\pi^n_j$ be the maximal path from  $(\fl{nx_j},\fl{nb_{k_j}})$
to $(\fl{nx_{j+1}},\fl{nb_{k_{j+1}+1}})$ in an environment with constant
weights $\omega_{i,j}=r_{\ell_j}^{-1}$ everywhere on the lattice,
with total weight $ T^n_j$.
 $ T^n_j\geq T(\pi^n_j)$, up to an error from the endpoints
of $\pi^n_j$.

Theorem 4.2 in \cite{sepp-large-deviations} gives
 a large deviation bound for  $ T^n_j$. Consider a constant 
rate $r$ environment and the maximal weight 
$T\bigl((\fl{nu_1},\fl{nv_1}), (\fl{nu_2},\fl{nv_2})\bigr)$
between two points $(u_1,v_1)$ and  $(u_2,v_2)$ such that 
their lattice versions can be connected by admissible paths for
all $n$. 
Then 
there exists a positive constant $C$ such that for $n$ large enough,
\be
\bP\Bigl\{ T\bigl((\fl{nu_1},\fl{nv_1}), (\fl{nu_2},\fl{nv_2})\bigr)
 > n r^{-1} 
\gamma( u_2-u_1 , v_2-v_1 ) + n\e  \Bigr\} < e^{-Cn}.
\label{ld1}
\ee 

There is a fixed finite collection out of which we pick 
the 
pairs $\{(x_j, b_{k_j}), (x_{j+1}, b_{k_{j+1}+1})\}$ that determine
the segments  $\tilde \pi^n_j$.
By \eqref{ld1} and  the Borel-Cantelli lemma, a.s.\ for  large enough $n$, 
\be
T^n_j \; \le\;  n r_{\ell_j}^{-1} 
\gamma( x_{j+1}-x_j , b_{k_{j+1}+1} - b_{k_j})
 + n\e \quad \text{for $0\le j<M_n$.}  
\label{ld3}
\ee  

Define $\delta_1>0$ by  $y+\delta_1=\sum_{j=0}^{M_n-1} (b_{k_{j+1}+1}-b_{k_j})$.
Since $y=\sum_{j=0}^{M_n-1} (b_{k_{j+1}}-b_{k_j})$ and by the choice
of the mesh of the partition $\{b_k\}$, we have 
$\delta_1\le  M_n\delta/Q\le \delta$.  
Think of $( x_{j+1}-x_j , b_{k_{j+1}+1} - b_{k_j})$,
$0\le j<M_n$,  as the successive
segments of a macroscopic  path from $(0,0)$  to $(x,y+\delta_1)$. 

For sufficiently large $n$ so that \eqref{ld3} is in effect,
\begin{align*} 
T^{n,\fl{nq}}(\fl{nx},\fl{ny}) &\le \sum_{j=1}^{M_n} T^n_j 
\le n \sum_{j=1}^{M_n}  r_{\ell_j}^{-1} 
\gamma( x_{j+1}-x_j , b_{k_{j+1}+1} - b_{k_j}) + nQ\e \\
&\le  n \Gamma^q(x,y+\delta_1)+ nQ\e  \\
&\le n\Gamma^q(x,y)+ nC\sqrt{\delta}  + nQ\e. 
\end{align*}
The last inequality came from Corollary \ref{deltapar2}. 
Let $\delta\rightarrow 0$. Since $\e$ was arbitrary the upper bound \eqref{ub7} holds. 
\end{proof}

 \begin{proof}[Proof of Theorem \ref{Tlimit}]  Fix $(x,y)$.  
For each  $\e  > 0$ we can find  lower semicontinuous step functions
$c_1$ and $c_2$ such that  $\norm{c_1-c_2}_\infty\le\e$ and on some
compact interval,  large enough to contain all the rates that can potentially
influence $\Gamma^q(x,y)$, 
   $c_1(x)\leq c(x)\leq c_2(x) $.    When the weights in \eqref{weights}  
   come  from speed function $c_i$  let us write 
  $T_i$ for last passage times and $\Gamma_i$ for their limits. 
An obvious
coupling using common exponential variables $\{\tau_{i,j}\}$ gives 
\[
T_1^{n,\fl{nq}}(\fl{nx},\fl{ny})\geq T^{n,\fl{nq}}(\fl{nx},\fl{ny})\geq T^{n,\fl{nq}}_2(\fl{nx},\fl{ny}).
\]   
Letting $\alpha>0$ denote a lower bound for $c(x)$ in the   compact interval relevant for $(x,y)$,
we have this bound for $\bx\in\cH(x,y)$: 
\begin{align*}
0\le \int_{0}^{1} \Bigl\{ \frac{\gamma(\textbf{x}'(s))}{c_1(x_1(s)-q)}
- \frac{\gamma(\textbf{x}'(s))}{c_2(x_1(s)-q)}\Bigr\} \,ds 
                             &\leq \e \int_{0}^{1}\frac{\gamma(\textbf{x}'(s))}{c^2_1(x_1(s)-q)}\,ds \\ 
                             &\leq \e  \alpha^{-2}\gamma(x,y).  
\end{align*} 
Therefore the limits also have the bound 
\begin{align*}
0\le \Gamma^q_1(x,y) - \Gamma^{q}_2(x,y) \leq C(x,y)\e. 
\end{align*}
From these approximations and the   limits for $T_i$  in Proposition \ref{simpleTlimit}
we can deduce Theorem \ref{Tlimit}. 
\end{proof}

\begin{proof}[Proof of Theorem \ref{two-phaselastpassagetime}]
We can construct the last passage times $G(x,y)$ of the
 corner growth model \eqref{basicref} with the same ingredients as the
wedge last passage times $T^{n,0}(x,y)$
of \eqref{lastpassagetime}, by taking $Y_{(i,j)}=\omega^{n,0}_{i-j,\,j}\tau^n_{i-j,\,j}$. 
Then  $T^{n,0}(x,y)=G(x+y,y)$ and  we can transfer the problem to the wedge.
The correct speed function to use  is now 
   $c(x) = c_1 \mathbf{1}\{ x < 0\}+c_2 \mathbf{1}\{ x\geq 0 \}$.
   In this case the limit in Theorem
\ref{Tlimit} can be solved explicitly with calculus.  We omit the details.  
 \end{proof}

\section{Hydrodynamic limit}
\label{hydrolimit}
In this section we sketch the proof of the main result Theorem \ref{hydro}. This argument is
from \cite{sepp99K, sepp01slow}.  

\subsection{Construction of the process and the variational coupling}
  For each $n \in \bN$ we construct a 
$\bZ$-valued \textsl{height process} $z^n(t)=(z^n_i(t): i \in \bZ)$.
 The height values obey the constraint 
\be
0\leq z^n_{i+1}(t)-z^n_i(t)\leq 1.
\label{constraint1}
\ee
Let $\{ \cD^n_i \}$ be a collection of mutually independent (in $i$ and $n$) Poisson processes with rates $c_i^n$ given by 
\be
c_i^n = c( n^{-1} i),
\label{discrete}
\ee
where $c(x)$ is the lower semicontinuous speed function. Dynamically, for each $n$ 
and  $i$, the height value $z^n_i$ is decreased by $1$ at  event times of $\cD^n_i$,
provided    the new configuration  does not violate \eqref{constraint1}. 
 
After we construct $z^{n}(t)$, we can define the exclusion process $\eta^{n}(t)$ by 
\be
\eta^n_i(t)= z^n_i(t)-z^n_{i-1}(t).
\label{exclusion-current}
\ee
A decrease in $z^n_i$ is associated with an exclusion  particle jump   from site $i$ to   $i+1$.  
Thus the $z^n$  process 
keeps track of the current of the $\eta^n $-process, precisely speaking  
 \be
J^n_i(t)=z^n_i(0)-z^n_i(t).
\label{zcurrent}
\ee 

Assume that the processes $z^n$ have  been constructed on a probability space that supports the initial configurations $z^n(0) = (z^n_i(0))$ and the Poisson processes $\{\cD^n_i\}$ that are independent of $(z^n_i(0))$.   Next we state the envelope property that is the key tool for the 
proof of the hydrodynamic limit.  Define a family of auxiliary height processes 
$\{\xi^{n,k}: n\in\bN,\,k\in\bZ\}$ that grow upward from  wedge-shaped initial conditions
\be
\xi_i^{n, k}(0) =\left\{
\begin{array}{ll}

\vspace{0.1 in}
0, & \textrm {if } i \geq 0  \\

-i, & \textrm{if } i < 0.
\end {array}
\right.
\label{xi-initial}
\ee
The dynamical rule for the $\xi^{n, k}$ process is that $\xi^{n,k}_i$ jumps up by 1
at the event times of $\cD^{n}_{i+k}$ provided the inequalities 
\be
\xi^{n,k}_i\leq\xi^{n,k}_{i-1} \quad \text{and} \quad \xi^{n,k}_i\leq\xi^{n,k}_{i+1}+1 
\label{xiineq}  
\ee
are not violated. In particular  $\xi^{n,k}_i$ attempts a jump at rate $c^n_{i+k}$.

\begin{lemma}[Envelope Property]
For each $n \in \bN$, for all $i\in \bZ$ and $t\geq 0$,
\be
z^n_i(t) = \sup_{k\in \bZ}  \{ z^n_k(0) - \xi^{n, k}_{i-k}(t)\} \quad a.s.
\label{envelope2}
\ee
\end{lemma}
Equation \eqref{envelope2}  holds by construction at time $t=0$, and it is
proved by induction on jumps. For  details see Lemma 4.2 in \cite{sepp99K}.

\subsection{The limit for $\xi$}
For $q,x \in \bR$, $t>0$ and for the speed function $c(x)$, define 
\be
g^{q}(x,t)= \inf\left\{ y: (x,y)\in \cW, \Gamma^q(x,y)\geq t  \right\}.
\label{gq}
\ee
$\Gamma^q(x,y)$ defined by \eqref{gammaq}  represents the macroscopic time it takes a $\xi$-type interface process 
to reach point $(x,y).$ 
The level curve of $\Gamma^q$ given by $g^q(\cdot,t)$ represents the limiting interface of a certain $\xi$-process,
as stated in the next proposition.
\begin{proposition}
For all $q,x \in \bR$ and $t>0$
\be
\dlim n^{-1}\xi^{n, \fl{nq}}_{\fl{nx}}(nt)=g^{-q}(x,t) \quad a.s.
\label{xilimit}
\ee 
\label{xilimit!}
\label{proxilimit}
\end{proposition}

Recall the lattice wedge  $\mathcal{L}$  defined by \eqref{defLL}.  
 For $(i,j)\in \cL\cup\partial\cL$, let 
\be
L^{n, k}(i,j)= \inf\{t\geq 0: \xi^{n,k}_i(t)\geq j\}
\ee
denote the time when $\xi^{n,k}_i$ reaches level $j$. The rules \eqref{xi-initial}--\eqref{xiineq} give the boundary conditions 
\be
L^{n,k}(i,j) = 0 \quad \text{for}\quad (i,j)\in \partial\cL
\label{Lboundary}
\ee and for $ (i,j) \in \cL$  the recurrence 
\be
L^{n,k}(i,j)= \max\{ L^{n,k}(i-1,j), L^{n,k}(i,j-1), L^{n,k}(i+1,j-1) \}+\beta^{n,k}_{i,j} 
\label{L}
\ee
where $\beta^{n,k}_{i,j}$ is an exponential waiting time, independent of everything else. It represents the time $\xi^{n,k}_i$ waits to jump, \textit{after} $\xi^{n,k}_i$ and its neighbors $\xi^{n,k}_{i-1}$, $\xi^{n,k}_{i+1}$ have reached positions that permit $\xi^{n,k}_i$ to jump
from $j-1$  to $j$. The dynamical rule that governs the jumps of $\xi^{n,k}_i$ implies that  $\beta^{n,k}_{i,j}$ has rate $c^n_{i+k}$.

Equations \eqref{lastpassagetime}, \eqref{boundaryoflastpassagetime}, \eqref{Lboundary},
and  \eqref{L}, together with the strong Markov property,  imply that 
\be
\{L^{n,k}(i,j):(i,j) \in \cL\cup\partial\cL \} \, \overset{\cD}= \, \{T^{n,-k}(i,j):(i,j) \in \cL\cup\partial\cL \}.
\label{L=-T}
\ee 
Consequently Theorem \ref{Tlimit}  gives 
the a.s.\ convergence $n^{-1} L^{n,\fl{nq}}(\fl{nx},\fl{ny})\to \Gamma^{-q}(x,y)$, 
and this passage time
limit gives limit \eqref{xilimit}.  

\begin{proof}[Proof of Theorem \ref{hydro}]
Given the initial configurations $\eta^n(0) = \left\{ \eta^n_i(0): i \in \bZ \right\}$ that appear in hypothesis \eqref{weaklaw}, define initial configurations $z^n(0) = \left\{ z^n_i(0): i \in \bZ \right\}$ so that $z^n_0(0)=0$ so that \eqref{exclusion-current} holds at time $t=0.$ Hypothesis \eqref{weaklaw} implies that
\be
\lim_{n\rightarrow \infty} n^{-1} z^n_{\fl{nq}} = v_0(q) \quad \textrm{a.s.} 
\label{weakz}
\ee  
for all $q \in \bR$, with $v_0$ defined by \eqref{vdef}.

Construct the height processes $z^n$ and define the exclusion processes $\eta^n$ by \eqref{exclusion-current}. Define $v(x,t)$ by \eqref{velocityversion1}. From 
\eqref{exclusion-current}--\eqref{zcurrent} we see that Theorem \ref{hydro}
 follows from proving that for all $x\in \bR, t\in \bR^+,$ 
\be
\lim_{n\rightarrow \infty}n^{-1}z^n_{\fl{nx}}(nt)= v(x,t) \quad \textrm{a.s.}
\label{z-v}
\ee  
Rewrite \eqref{envelope2} with the correct scaling:
\be
n^{-1}z^n_{\fl{nx}}(nt)= \sup_{q \in \bR}\left\{ n^{-1}z^n_{\fl{nq}}(0) - n^{-1}\xi^{\fl{nq}}_{\fl{nx}-\fl{nq}}(nt) \right\}.
\label{scaledz}
\ee

The proof of \eqref{z-v} is now to show that the right-hand side of \eqref{scaledz} converges to the right-hand side of \eqref{velocityversion1}.
 
From \eqref{weakz}, \eqref{scaledz} and \eqref{xilimit} we can prove that a.s.
\be
\lim_{n\rightarrow \infty}n^{-1}z^n_{\fl{nx}}(nt)= \sup_{q \in \bR}\left\{ v_0(q)-g^{-q}(x-q, t) \right\}\equiv \tilde{v}(x,t).
\label{tildev}
\ee 
The argument is the same as the one from equations (6.4)--(6.15) in \cite{sepp99K} so we will not repeat it here. 

Using \eqref{gammaq} and \eqref{gq} we can rewrite $\tilde{v}(x,t)$ as
\be
\tilde{v}(x,t)= \sup_{q,y \in \bR}\Big\{ v_0(q)-y: \exists \bx \in \cH(x-q,y) \textrm{ such that } \int_0^1\frac{\gamma(\bx'(s))}{c(x_1(s) + q)}\,ds \geq t \Big\}.
\label{massagedvtilde}
\ee  
The final step is to prove $v(x,t) = \tilde{v}(x,t).$ The argument is identical to the one used to prove Proposition 4.3 in \cite{sepp01slow} so we omit it.   With this we can consider 
Theorem \ref{hydro} proved.  
\end{proof}
  
\section{Density profiles in two-phase TASEP}
\label{density}
This section proves Corollary \ref{densityprofiles}: assuming  $c(x) =(1- H(x))c_1 + H(x)c_2$, $c_1 \geq c_2$ and $\r_0(x) \equiv \r \in (0,1)$,  we   use  variational formula \eqref{velocityversion1}
to obtain   explicit hydrodynamic limits. 
\begin{remark}
In light of Corollary \ref{pdecor}, one can (instead of doing the following computations) guess the candidate solution for the scalar conservation law  \eqref{Breakscl} and then check that it verifies the entropy conditions \eqref{Eb1} - \eqref{Eb3}. The following computations do not require any knowledge of p.d.e.\ theory or familiarity with interface problems so we present them independently in this section.
\end{remark}
 Let 
\be
C^0(x,t,q) = \left\{w\in C( [0,t], \bR):  w \text{ piecewise linear, } w(0)=q, w(t)=x \right\}.  
\ee
To   optimize in \eqref{velocityversion1} we use a couple
different approaches for different cases. We outline this and omit the details.  

One approach is to separate the choice of the starting point $q$ of the path. 
By setting  
\be
I(x,t,q) = \inf_{w \in C^0(x,t,q)}\bigg\{ \int_{0}^{t}c(w(s))g\left( \dfrac{w'(s)}{c(w(s))} \right)\,ds \bigg\}
\label{Ixtq}
\ee
  \eqref{velocityversion1} becomes\be
v(x,t) = \sup_{q \in \bR}\left\{ v_0(q) - I(x,t,q) \right\}.
\label{spoptimize}
\ee
We  distinguish four cases according to the signs of $x,q$.  Set 
\be
R_{+}(x,t) = \sup_{q > 0}\left\{ v_0(q) - I(x,t,q) \right\},\textrm{ if } x>0,\quad
\ee
\be
L_{-}(x,t) = \sup_{q < 0}\left\{ v_0(q) - I(x,t,q) \right\},\textrm{ if } x<0.
\ee 
These functions are going to be used in Cases 1 and 2 below ($qx \geq 0$) where we can 
compute $I(x,t,q)$ directly. 

However, there are values  $(x,t, q)$ for which the $q$-derivative of the 
expression in   braces in \eqref{spoptimize} is a rational function with a quartic polynomial in the numerator. While an explicit formula for   roots of a quartic  exists, the solution is not attractive and
it is not clear how to pick the right root. Instead we turn the problem into a two-dimensional maximization problem. 

If $qx < 0$    the optimizing path $w$ crosses the origin:
$w(u)= 0$ for some $u$.  It turns out convenient  
to find the optimal $q$ for each crossing time $u$. For Case 3 ($q<0,x>0$) set
\be
\Phi(u, q) = q\r - c_1ug\left(\dfrac{-q}{uc_1}\right) - c_2(t-u)g\left( \dfrac{x}{(t-u)c_2} \right)
\ee 
and 
\be
L_{+}(x,t) = \displaystyle\sup_{q < 0, u \in[0,t]}\Phi(u,q).  
\ee
For Case 4 ($q>0,x<0$)   the obvious modifications are 
\be
\Psi(u, q) = q\r - c_2ug\left(\dfrac{-q}{uc_2}\right) - c_1(t-u)g\left( \dfrac{x}{(t-u)c_1} \right)
\ee 
and
\be
R_{-}(x,t) = \displaystyle\sup_{q > 0, u \in[0,t]}\Psi(u,q).  
\notag
\ee
Rewrite \eqref{spoptimize} using   functions $R_{\pm}$, $L_{\pm}$: 
\be
v(x,t)= \max\{ R_{+}(x,t), L_{+}(x,t)\}\mathbf{1}\{x\geq 0\} + \max\{ R_{-}(x,t), L_{-}(x,t)\}\mathbf{1}\{x < 0\}.
\label{R-Lvelocity} 
\ee

\begin{proof}[Proof of Corollary \ref{densityprofiles}\\]

We  compute the functions  $R_{\pm}$, $L_{\pm}$.
 The density profiles $\r(x,t)$ are given then by the $x$-derivative of $v(x,t)$.  
\bigskip

\noindent\textsl{\textbf{Case 1}: $x\geq 0$, $q \geq 0$.}
Since $c_2 \leq c_1$, the minimizing $w$ of $I(x,t,q)$ is the straight line connecting $(0,q)$ to $(t,x).$ In particular, 
\be I(x,t,q)=c_2tg\left(\frac{x-q}{tc_2}\right).\ee Then the resulting $R_{+}(x,t)$ is given by
\be
R_{+}(x,t) =\begin{cases}

-tc_2g(\frac{x}{tc_2} ) &\textrm{if } \r\leq \frac{1}{2},\quad x< tc_2(1-2\r)\\[4pt]

\r x -tc_2\r(1-\r), & \textrm{if } \r\leq \frac{1}{2},\quad x\geq tc_2(1-2\r)\\[4pt]

\r x -tc_2\r(1-\r), & \textrm{if } \r > \frac{1}{2}, 
\end{cases}
\ee 

\noindent\textsl{\textbf{Case 2}: $x\leq 0$, $q \leq 0$.}
The minimizing path $w$ can either be a straight line from $(0,q)$ to $(t,x)$ or a piecewise linear path such that the set $\{t: w(t) = 0\}$ has positive Lebesgue measure. This last statement just says that the path might want to take advantage of the low rate at $x=0$. We leave the calculus details to the reader and record the resulting minimum value of $I(x,t,q)$.  Set   $B = \sqrt{c_1(c_1-c_2)}$.  
\be
I(x,t,q) =
\begin{cases}
\vspace{0.07in}
\frac{-qc_1}{4B}\left(1-\frac{B}{c_1} \right)^2 + \left( t-\frac{|x|-q}{B} \right)\frac{c_2}{4} -\frac{xc_1}{4B}\left(1+\frac{B}{c_1} \right)^2,\\

 \quad\quad\quad\quad\quad\quad\quad\quad  \text{when } -(\sqrt{Bt}-\sqrt{|x|})^2 \leq q,\, -Bt\leq x< 0\quad\\ 

c_1tg\left(\frac{x-q}{c_1t}\right) \quad\quad\quad\text{otherwise }
\end{cases}
\ee 
 The corresponding function $L_{-}(x,t)$ is given by\be
L_{-}(x,t) =
\begin{cases}

\vspace{0.07 in}
\r x -tc_1\r(1-\r), &0< \r < \r^*, x\in \bR\\

\vspace{0.07 in}
\r x -tc_1\r(1-\r), &\r^*\leq\r\leq \frac{1}{2},\quad x\leq -tc_1(\r-\r^*)\\ 

\vspace{0.07 in}
-\left( t+\frac{x}{B} \right)\frac{c_2}{4} 
+  \frac{xc_1}{4B}\left(1+\frac{B}{c_1} \right)^2, & \r^*\leq\r\leq \frac{1}{2},\quad x > -tc_1(\r-\r^*)\\  

\vspace{0.07 in}
\r x -tc_1\r(1-\r) ,& \frac{1}{2} < \r \leq 1 - \r^*,\quad x< -tc_1(\r-\r^*) \\ 

\vspace{0.07 in}
-\left( t+\frac{x}{B} \right)\frac{c_2}{4} 
+  \frac{xc_1}{4B}\left(1+\frac{B}{c_1} \right)^2, & \frac{1}{2} < \r \leq 1 - \r^*,\quad -tc_1(\r-\r^*)\leq x \\

\vspace{0.07 in}
-\left( t+\frac{x}{B} \right)\frac{c_2}{4} 
+  \frac{xc_1}{4B}\left(1+\frac{B}{c_1} \right)^2, &1 - \r^* < \r<1,\quad -Bt\leq x \\

\vspace{0.07 in}
-tc_1g\big( \frac{x}{tc_1} \big), & 1-\r^* < \r<1, \quad  -c_1t(2\r-1) \leq x < -Bt \\

\vspace{0.07 in}
\r x - c_1t\r(1-\r), & 1-\r^* < \r<1, \quad  x < -c_1t(2\r-1) \\

\end{cases}
\ee 

\medskip

\noindent\textsl{\textbf{Case 3}: $x > 0$, $q \leq 0.$}
Abbreviate $ D = c_2^2 - 4c_1c_2\r(1-\r)$. First compute the $q$-derivative\be
\Phi_q(u,q) =
\begin{cases}

\vspace{0.07 in}
\r - \frac{1}{2} - \frac{q}{2uc_1}, & -uc_1 \leq q < 0 \\ 

\r & q < -uc_1.
\end{cases}
\ee 
If $\r\geq 1/2$ then $\Phi_q$ is positive and the maximum value is when $q=0$ so we are reduced to Case $1$. If $\r< 1/2$ the maximizing $\displaystyle q = uc_1\left(2\r - 1\right).$
Then 
\be F(u) = \Phi\bigl(u, 2uc_1(\r - \tfrac{1}{2})\bigr) = -uc_1\r(1-\r) -c_2(t-u)g\left( \dfrac{x}{(t-u)c_2} \right),\notag\ee with $u$-derivative 
\be
\frac{dF}{du} =  -c_1\r(1-\r)+ \frac{c_2}{4}\left(1-\frac{x^2}{(c_2(t-u))^2}\right).
\notag\ee 
Again we need to split two cases. 
If $\r < \r^*$ (equivalently $D>0$) and $x\leq t\sqrt{D}$, the maximizing $u = t- x/\sqrt{D},$ otherwise $u=0$. If $\r^* \leq \r < \frac{1}{2}$ the derivative is negative so the maximizing $u$ is still $u=0.$  Together,
\be
L_{+}(x,t) =
\begin{cases}

\vspace{0.07 in}
-tc_1\r(1-\r)+x\Big(\frac{1}{2}-\frac{\sqrt{D}}{2c_2} \Big), & \r<\r^*, x\leq t\sqrt{D}\\

\vspace{0.07 in}
-c_2tg\big(\frac{x}{tc_2}\big) , &  \r<\r^*, x\geq t\sqrt{D}\\ 

-c_2tg\big(\frac{x}{tc_2}\big) , & \r^*\leq\r \leq 1. 

\end{cases}
\ee 

\medskip

\noindent\textsl{\textbf{Case 4}: $x \leq 0, q \geq 0.$}
We treat this case in exactly the same way as Case 3, so we omit the details. Here we need the quantity $D_1=(c_1)^2 - 4c_1c_2\r(1-\r)$ and we compute 
\be
R_{-}(x,t) =
\begin{cases}

\vspace{0.07 in}
-tc_1g\big(\frac{x}{tc_1}\big), & \r\leq\frac{1}{2}\\

\vspace{0.07 in}
-tc_2\r(1-\r)+x\Big( \frac{1}{2}+\frac{\sqrt{D_1}}{2c_1} \Big) , & \frac{1}{2}<\r,- t\sqrt{D_1}\leq x\\

-tc_1g\big(\frac{x}{tc_1}\big) , & \frac{1}{2} \leq \r \leq 1, x<- t\sqrt{D_1}\\ 
\end{cases}
\ee
Now compute $v(x,t)$ from \eqref{R-Lvelocity}. We leave the remaining details to the reader.
\end{proof}

\section{Entropy solutions of the discontinuous conservation law}
\label{pdes}

For this section, $c(x)=(1-H(x))c_1 + H(x)c_2$, $h(\r)=\r(1-\r)$ and set $F(x,\r) = c(x)h(\r)$ for the flux function of the scalar conservation law \eqref{scl} and $\widetilde{F}(x,\r) =c(x)f(\r)$ for the flux function of the particle system, where $f$ is given by \eqref{fnot}. (The difference between
$F$ and $\wt F$ is that the latter is $-\infty$ outside $0\le \rho\le 1$.) 
 
In \cite{adim-gowd-03} the authors prove that there exists a solution to the corresponding Hamilton-Jacobi equation
\be
\left\{
\begin{array}{lll}

\vspace{0.1 in}
V_t+c_1h(V_x) =0, & \textrm{if } x <0, t>0 \\

\vspace{0.1 in}
V_t+c_2h(V_x) =0, & \textrm {if } x > 0, t>0  \\

V(x,0)=V_0(x)
\end {array}
\right.
\label{HJ}
\ee 
such that $V_x$ solves the scalar conservation law \eqref{scl} with flux function $F(x,\rho)$ and $V_x$ satisfies the entropy assumptions $(E_i), (E_b)$. $V(x,t)$ is given by 
\be
V(x,t)=\sup_{w(\cdot)}\left\{  V_0(w(0)) + \int_{0}^{t} (c(w(s))h)^*(w'(s)) \,ds \right\},
\label{HJsol}
\ee 
where the supremum is taken over  piecewise linear paths $w\in C([0,t],\bR)$ that satisfy $w(t)=x.$

To apply the results of \cite{adim-gowd-03} to the profile
$\rho(x,t)$ coming from our hydrodynamic limit, we only need
to show that the variational descriptions match, in other words
that  we can replace $F$ with $\widetilde{F}$ and 
the solution is still the same.

\begin{proof}[Proof of Theorem \ref{pdecor}]
Convex duality gives  $\left(c(x)f\right)^*(y) = 
c(x)f^*\left( y/c(x)\right)$ and 
so we can rewrite  \eqref{velocityversion1} as 
\be
v(x,t)=\sup_{w(\cdot)}\left\{  v_0(w(0)) + \int_{0}^{t} (c(w(s))f)^*(w'(s)) \,ds \right\}.
\label{TNsol}
\ee 
Observe that for all $y \in \bR$
\be
(c(x)f)^*(y) \geq  (c(x)h)^*(y),
\label{flux-comparison}
\ee
with equality if and only if $y \in [-c_1,c_2]$
Since the supremum in \eqref{HJsol} and \eqref{TNsol} is taken over the same set of paths, \eqref{flux-comparison} implies that
\be
V(x,t) \leq v(x,t).
\label{v-comparison}
\ee 
The proof of the theorem is now reduced to proving that the supremum in \eqref{TNsol} is achieved when $w'(s)c(w(s))^{-1} \in [-1,1],$ giving $V(x,t)=v(x,t)$.

To this end we rewrite $v(x,t)$ once more,
this time as
\[
v(x,t)= \max\{ R_{+}(x,t), L_{+}(x,t)\}\mathbf{1}\{x\geq 0\} + \max\{ R_{-}(x,t), L_{-}(x,t)\}\mathbf{1}\{x < 0\} 
\]
where the functions $R_{\pm}$, $L_{\pm}$ 
(as in the proof of Corollary \ref{densityprofiles})
are defined by 
\be
R_{+}(x,t) = \sup_{q > 0}\left\{ v_0(q) - I(x,t,q) \right\},\textrm{ if } x>0,\quad
\ee
\be
L_{-}(x,t) = \sup_{q < 0}\left\{ v_0(q) - I(x,t,q) \right\},\textrm{ if } x<0,
\ee 
where $I(x,t,q)$ is as in  \eqref{Ixtq}, 
and 
\be
L_{+}(x,t) = \displaystyle\sup_{q < 0, u \in[0,t]} \bigg\{ v_0(q) - c_1u g\Big(\frac{-q}{uc_1}\Big) - c_2(t-u)g\Big(\frac{x}{(t-u)c_2} \Big)\bigg\} \quad \textrm{if } x\geq 0,
\label{finalL}
\ee
and 
\be
R_{-}(x,t) = \displaystyle\sup_{q > 0, u \in[0,t]}\bigg\{ v_0(q) - c_2ug\left(\dfrac{-q}{uc_2}\right) - c_1(t-u)g\left( \dfrac{x}{(t-u)c_1} \right) \bigg\}, \quad x\le 0.  
\ee

 It suffices to show that the suprema
 that define  $R_{\pm}$, $L_{\pm}$
are achieved when \be w'(s)c(w(s))^{-1} \in [-1,1].
\label{slopecondition}
\ee 
%for arbitrary initial conditions $\r_0(x) \in [0,1].$ 
We show this for $L_{+}$. The remaining cases are similar. In \eqref{finalL},
 as before, $u$ is the time for which $w(u) = 0.$ 
Let $\Phi(u,q)$ denote the expression in braces in \eqref{finalL}
with  $q$-derivative
\be
\Phi_q(u,q) =
\begin{cases}

\vspace{0.07 in}
\r_0(q) - \frac{1}{2} - \frac{q}{2uc_1}, &-uc_1 \leq q < 0 \\ 

\r_0(q), & q < -uc_1.
\end {cases}
\label{qder}
\ee 

Observe that if $\Phi_q(u,q) = 0$ for some $q^* = q^*(u)$ then also $q^*$ maximizes $\Phi$.
Otherwise the maximum is achieved at $0$ and we are reduced to a different case. Assume that $q^*$ exists. Then by \eqref{qder}
\be
\frac{-q^*}{u} = (1- 2\r_0(q^*))c_1 < c_1.
\label{slope1}
\ee  
Therefore, the slope of the first segment of the maximizing path $w$ satisfies \eqref{slopecondition}. 

The slope of the second segment is $x(t-u)^{-1}.$ Assume  that the piecewise linear path $w$ defined by $u$ and $q^*$ is the one that achieves the supremum. Also assume $u > t-xc_2^{-1}.$ Consider the path $\tilde{w}$ with $\tilde{w}(0) = q^*$ and $\tilde{w}(t-xc_2^{-1})= 0$. Since $g$ is decreasing, we only increase the value of $\Phi$. Hence the supremum that gives $L_{+}$ cannot be achieved on $w$ and this gives the desired contradiction.
\end{proof}
\medskip

\textbf{Acknowledgments.} We would like to thank an anonymous referee for thoroughly reading the original manuscript and for offering comments that improved the presentation of this paper.   

%\bibliography{refs,growthrefs,pderefs} %miktek 2.7
%\bibliographystyle{plain}

\end{document}